\documentclass[10pt]{article} 
\usepackage{geometry}                       
\usepackage{graphicx}
\usepackage[table]{xcolor}             
\usepackage{amssymb}
\usepackage{amsmath}
\usepackage{mathtools} 
\usepackage{svg}
\usepackage[utf8]{inputenc}
\usepackage{amsthm}
\usepackage{bm}
\usepackage{xcolor}

\usepackage{booktabs}
\usepackage{amsthm}

\newtheorem{exmp}{Example}

\newtheorem{rmrk}{Remark}

\newtheorem{algo}{Algorithm}
\newtheorem*{wh}{Working Hypothesis}

\usepackage{MnSymbol}
\usepackage{algorithm}
\usepackage[noend]{algpseudocode}
\usepackage{cite}
\usepackage{booktabs}
\usepackage{listings}
\usepackage[]{hyperref}
\hypersetup{
    colorlinks = true,
    linkbordercolor = {white},
    linkcolor=black,
    citecolor=black,
    urlcolor=black
}
\makeatletter
\newcommand{\doubletilde}[1]{{%
  \mathpalette\double@tilde{#1}%
}}
\newcommand{\double@tilde}[2]{%
  \sbox\z@{$\m@th#1\tilde{#2}$}%
  \ht\z@=.9\ht\z@
  \tilde{\box\z@}%
}
\makeatother
\usepackage{amsthm}
\usepackage[export]{adjustbox}
\usepackage{mathdots}

\definecolor{jlbase}{HTML}{444444}            
\definecolor{jlkeyword}{HTML}{444444}         
\definecolor{jlliteral}{HTML}{78A960}         
\definecolor{jlbuiltin}{HTML}{397300}         
\definecolor{jlcomment}{HTML}{888888}         
\definecolor{jlstring}{HTML}{880000}          
\definecolor{jlbackground}{HTML}{F5F5F5} 
\lstdefinelanguage{Julia}%
  {morekeywords=[1]{abstract,break,case,catch,const,continue,do,else,elseif,%
      end,export,false,for,function,immutable,import,importall,if,in,%
      macro,module,otherwise,quote,return,switch,true,try,type,typealias,%
      using,while},%
   morekeywords=[2]{LinearAlgebra,GenericLinearAlgebra},%
    morekeywords=[3]{diagm,zeros,ones,cos,range,LinRange,sort,length,setprecision,big,real,convert,isreal,eigvals},
        morekeywords=[4]{Float64,BigFloat,Bool,Integer,Array},
   sensitive=true,%
   alsoother={$},%
   morecomment=[l]\#,%
   morecomment=[n]{\#=}{=\#},%
   morestring=[s]{"}{"},%
   morestring=[m]{'}{'},%
}[keywords,comments,strings]%

\begin{document}
\title{A Matrix-Less Method to Approximate the Spectrum and the Spectral Function of Toeplitz Matrices with Real Eigenvalues}
\author{Sven-Erik Ekström\\ {\small\texttt{see@2pi.se}}}
\date{%
    Athens University of Economics and Business\\[2ex]%
    \today
}
\maketitle
\begin{abstract}
It is known that the generating function $f$ of a sequence of Toeplitz matrices $\{T_n(f)\}_n$ may not describe the asymptotic distribution of the eigenvalues of $T_n(f)$ if $f$ is not real. In this paper, we assume as a working hypothesis that, if the eigenvalues of $T_n(f)$ are real for all $n$, then they admit an asymptotic expansion of the same type as considered in previous works \cite{ekstrom171,ahmad171,ekstrom183,ekstrom185}, where the first function $g$ appearing in this expansion is real and describes the asymptotic distribution of the eigenvalues of $T_n(f)$. After validating this working hypothesis through a number of numerical experiments, drawing inspiration from \cite{ekstrom183}, we propose a matrix-less algorithm in order to approximate the eigenvalue distribution function $g$. 
The proposed algorithm is tested on a wide range of numerical examples; in some cases, we are even able to find the analytical expression of $g$. Future research directions are outlined at the end of the paper.
\end{abstract}

\section{Introduction}
\label{sec:introduction}
Given a function $f\in L^1([-\pi,\pi])$, the $n\times n$ Toeplitz matrix generated by $f$ is defined as 
\begin{align}
T_n(f)=\left[
\begin{array}{ccccccccc}
\hat{f}_0&\hat{f}_{-1}&\cdots&\hat{f}_{n-1}\vphantom{\vdots}\\
\hat{f}_1&\ddots&\ddots&\vdots\\
\vdots&\ddots&\ddots&\hat{f}_{-1}\\
\hat{f}_{n-1}&\cdots&\hat{f}_1&\hat{f}_0\vphantom{\vdots}\\
\\[-0.6em]
\end{array}
\right],\nonumber
\end{align}%
where the numbers $\hat{f}_k$ are the Fourier coefficients of $f$, that is,
\begin{align}
\hat{f}_k=\frac{1}{2\pi}\int_{-\pi}^{\pi}f(\theta)e^{-\mathbf{i}k\theta}\mathrm{d}\theta,\qquad f(\theta)=\sum_{k=-\infty}^{\infty}\hat{f}_ke^{\mathbf{i}k\theta}.\label{eq:introduction:fourier}
\end{align}
It is known that the generating function $f$, also known as the symbol of $\{T_n(f)\}_n$,
describes the asymptotic distribution of the singular values of $T_n(f)$; if $f$ is real or if $f\in L^\infty([-\pi,\pi])$ and its essential range has empty interior and does not disconnect the complex plane, then $f$ also describes the asymptotic distribution of the eigenvalues of $f$; see \cite{bottcher991,garoni171,tilli991} for details and~\cite[Section~3.1]{garoni171} for the notion of asymptotic singular value and eigenvalue distribution of a sequence of matrices. We will write $\{T_n(f)\}_n\sim_\sigma f$ to indicate that $\{T_n(f)\}_n$ has an asymptotic singular value distribution described by $f$ and $\{T_n(f)\}_n\sim_\lambda f$ to indicate that $\{T_n(f)\}_n$ has an asymptotic eigenvalue distribution described by $f$. The cases of interest in this paper are those in which $\{T_n(f)\}_n\not\sim_\lambda f$ and the eigenvalues of $T_n(f)$ are real for all $n$. We believe that in these cases there exist a real function $g$ such that $\{T_n(f)\}_n\sim_\lambda g$ and the eigenvalues of $T_n(f)$ admit an asymptotic expansion of the same type as considered in previous works \cite{ekstrom171,ahmad171,ekstrom183,ekstrom185}. We therefore formulate the following working hypothesis. 

\begin{wh}
Suppose that the eigenvalues of $T_n(f)$ are real for all $n$. 
Then, for every integer $\alpha\ge0$, every $n$ and every $j=1,\ldots,n$, the following asymptotic expansion holds:
\begin{align}\label{eq:introduction:hoapp}
\lambda_j(T_n(f))&=g(\theta_{j,n})+\sum_{k=1}^\alpha c_k(\theta_{j,n})h^k+E_{j,n,\alpha},\nonumber\\
&=\sum_{k=0}^\alpha c_k(\theta_{j,n})h^k+E_{j,n,\alpha},
\end{align}
where:
\begin{itemize}
	\item the eigenvalues of $T_n(f)$ are arranged in non-decreasing order, $\lambda_1(T_n(f))\le\ldots\le\lambda_n(T_n(f))$;
	\item $\{g\coloneqq c_0,c_1,c_2,c_3,\ldots\}$ is a sequence of functions from $(0,\pi)$ to $\mathbb R$ which depends only on $f$;
	\item $h=\frac{1}{n+1}$ and $\theta_{j,n}=\frac{j\pi}{n+1}=j\pi h$;
	\item $E_{j,n,\alpha}=O(h^{\alpha+1})$ is the remainder (the error), which satisfies the inequality $|E_{j,n,\alpha}|\le C_\alpha h^{\alpha+1}$ for some constant $C_\alpha$ depending only on $\alpha,f$.
\end{itemize}
\end{wh}

\begin{rmrk}
In the working hypothesis, we arrange the eigenvalues of $T_n(f)$ in non-decreasing order, however, using a non-increasing order would result in another function $g$. The case where the eigenvalues of $T_n(f)$ can be described by a complex-valued or non-monotone function $g$ is out of the scope of this article and warrants further research.
\end{rmrk}

\section{Motivation and illustrative examples}
\label{sec:motivation}
In this section we present four examples in support of our working hypothesis.
We also discuss the fact that standard double precision eigenvalue solvers (such as LAPACK, \texttt{eig} in \textsc{Matlab}, and \texttt{eigvals} in \textsc{Julia}) fail to give accurate eigenvalues of certain matrices $T_n(f)$; see, e.g., \cite{beam931,trefethen051}. High-precision computations, by using packages such as \textsc{GenericLinearAlgebra.jl}~\cite{noack191} in \textsc{Julia} can compute the true eigenvalues, but they are very expensive from the computational point of view. 
Therefore, approximating $g$ on the grid $\theta_{j,n}$ and using matrix-less methods~\cite{ekstrom183} to compute the spectrum of $T_n(f)$ can be computationally very advantageous. Also, the presented approaches can be a valuable tool for the analysis of the spectra of non-normal Toplitz matrices having real eigenvalues. 

Here is a short description of the four examples we are going to consider.
In what follows, we denote by $\xi_{j,n}$ a ``perfect'' sampling grid, typically not equispaced, such that $\lambda_j(T_n(f))=g(\xi_{j,n})$ for $j=1,\ldots,n$; such grids are discussed in~\cite{ekstrom191}.
\begin{itemize}
\item Example~\ref{exmp:1}: $T_n(f)$ is non-symmetric tridiagonal, $g$ is known, and the eigenvalues $\lambda_j(T_n(f))=g(\theta_{j,n})$ are known explicitly;
\item Example~\ref{exmp:2}: $T_n(f)$ is symmetric pentadiagonal, $g=f$, and the eigenvalues $\lambda_j(T_n(f))=g(\xi_{j,n})$ are not known explicitly;
\item Example~\ref{exmp:3}: $T_n(f)$ is non-symmetric, $g$ is known, and the eigenvalues $\lambda_j(T_n(f))=g(\xi_{j,n})$ are not known explicitly;
\item Example~\ref{exmp:4}: $T_n(f)$ is non-symmetric, $g$ is not known, and the eigenvalues $\lambda_j(T_n(f))=g(\xi_{j,n})$ are not known explicitly.
\end{itemize}

\begin{exmp}
\label{exmp:1}
Consider the symbol
\begin{align}
f(\theta)&=\hat{f}_{1}e^{\mathbf{i}\theta}+\hat{f}_0+\hat{f}_{-1}e^{-\mathbf{i}\theta}.
\label{eq:exmp:1:ftridiag}
\end{align}
The matrix $T_n(f)$ is tridiagonal, and there exist a function
\begin{align}
g(\theta)=\hat{f}_0+2\sqrt{\hat{f}_{1}}\sqrt{\hat{f}_{-1}}\cos(\theta),
\label{eq:exmp:1:gtridiag}
\end{align}
such that $T_n(f)\sim T_n(g)$, that is, they are similar and hence have the same eigenvalues. The eigenvalues are given explicitly by
\begin{align}
\lambda_j(T_n(f))=g(\theta_{j,n}),
\label{eq:exmp:1:eigexacttridiag}
\end{align}
where $\theta_{j,n}$ is defined in the working hypothesis.  
Now, choose the Fourier coefficients $\hat{f}_{1}=-1$, $\hat{f}_0=2$, and $\hat{f}_{-1}=-2$. In this case, we have
\begin{align}
f(\theta)&=-e^{\mathbf{i}\theta}+2-2e^{-\mathbf{i}\theta},\nonumber\\
g(\theta)&=2-2\sqrt{2}\cos(\theta),
\label{eq:exmp:1:gexmptridiag}
\end{align}
and the spectrum of $T_n(f)$ is real, even though the symbol $f$ is complex-valued.
The Toeplitz matrices generated by $f$ and $g$ are given by
\begin{align}
T_n(f)=\left[
\begin{array}{rrrrr}
2&-2\\
-1&2&-2\\
&\ddots&\ddots&\ddots\\
&&\ddots&\ddots&-2\\
&&&-1&2
\end{array}
\right],
\quad
T_n(g)=\left[
\begin{array}{rrrrr}
2&-\sqrt{2}\\
-\sqrt{2}&2&-\sqrt{2}\\
&\ddots&\ddots&\ddots\\
&&\ddots&\ddots&-\sqrt{2}\\
&&&-\sqrt{2}&2
\end{array}
\right],\nonumber
\end{align}
We also note that $T_n(g)$ is a symmetrized version of $T_n(f)$, in the sense that there exists a decomposition $T_n(g)=DT_n(f)D^{-1}$ where $D$ is a diagonal matrix with elements $(D)_{i,i}=\gamma^{i-1}$, and $\gamma=\sqrt{\hat{f}_{-1}}/\sqrt{\hat{f}_1}$; see \cite{parter621}.

\begin{figure}[!ht] 
\centering
\includegraphics[width=0.472\textwidth,valign=t]{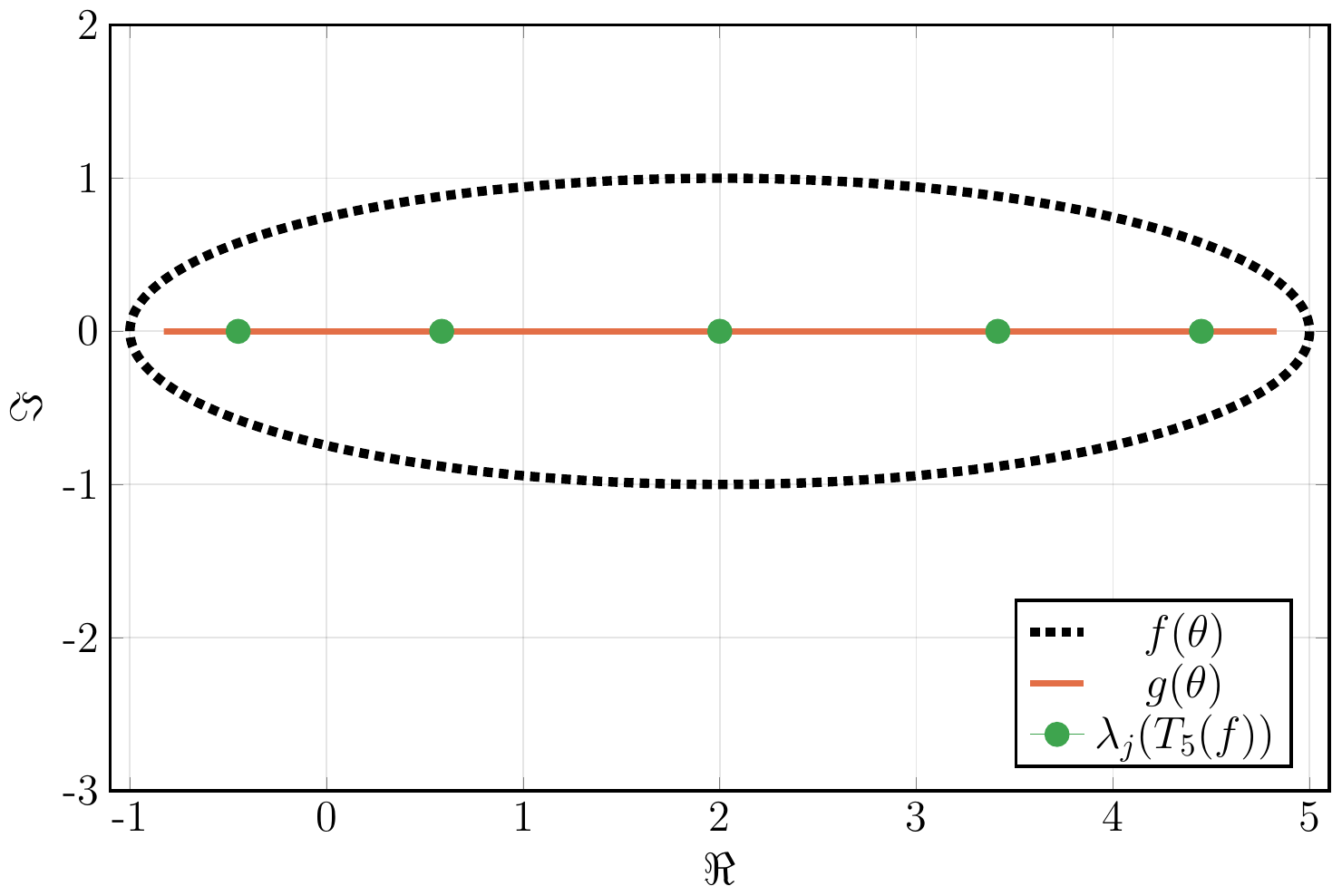}
\includegraphics[width=0.48\textwidth,valign=t]{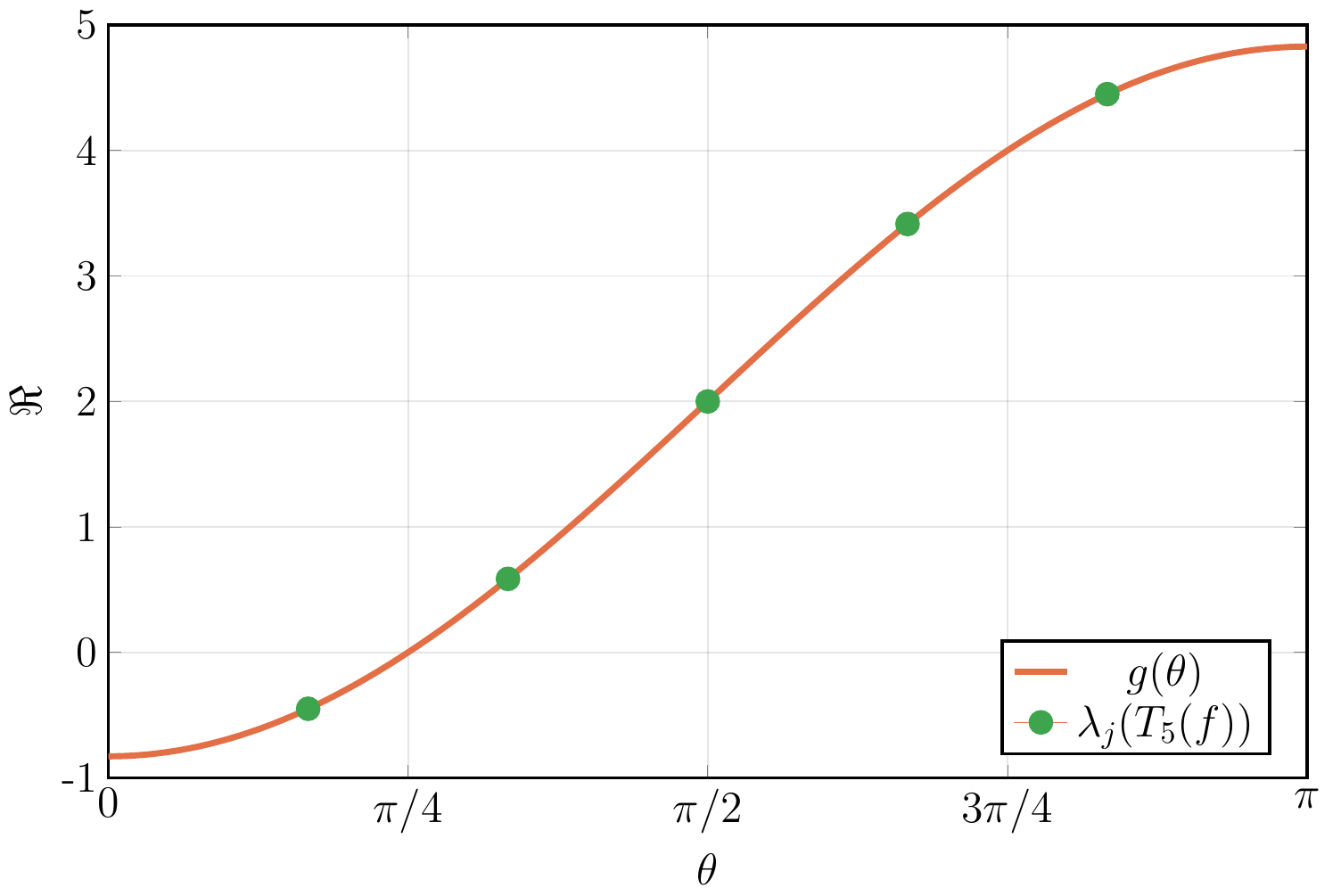}
\caption{[Example~\ref{exmp:1}: Symbol $f(\theta)=-e^{\mathbf{i}\theta}+2-2e^{-\mathbf{i}\theta}$] Left: Representations of $f(\theta)$ (dashed black line), and $g(\theta)=2-2\sqrt{2}\cos(\theta)$ (red line), and $\lambda_j(T_n(f))=\lambda_j(T_n(g))$ for $n=5$ (green dots). Right: Representations of $g$ and $\lambda_j(T_5(f))=\lambda_j(T_5(g))=g(\theta_{j,n})$.}
\label{fig:exmp:1:symbols}
\end{figure}
In the left panel of Figure~\ref{fig:exmp:1:symbols} we represent the function $f$ (dashed black line), and $g$ (red line), and the eigenvalues $\lambda_j(T_n(f))=\lambda_j(T_n(g))$ (green dots) for $n=5$. 
In the right panel of Figure~\ref{fig:exmp:1:symbols}  we show the function $g$ (red line) on the interval $[0,\pi]$ only (since it is even on $[-\pi,\pi]$) and the eigenvalues $\lambda_j(T_5(f))=\lambda_j(T_5(g))=g(\theta_{j,5})$ (green dots).
\begin{figure}[!ht]
\centering
\includegraphics[width=0.48\textwidth]{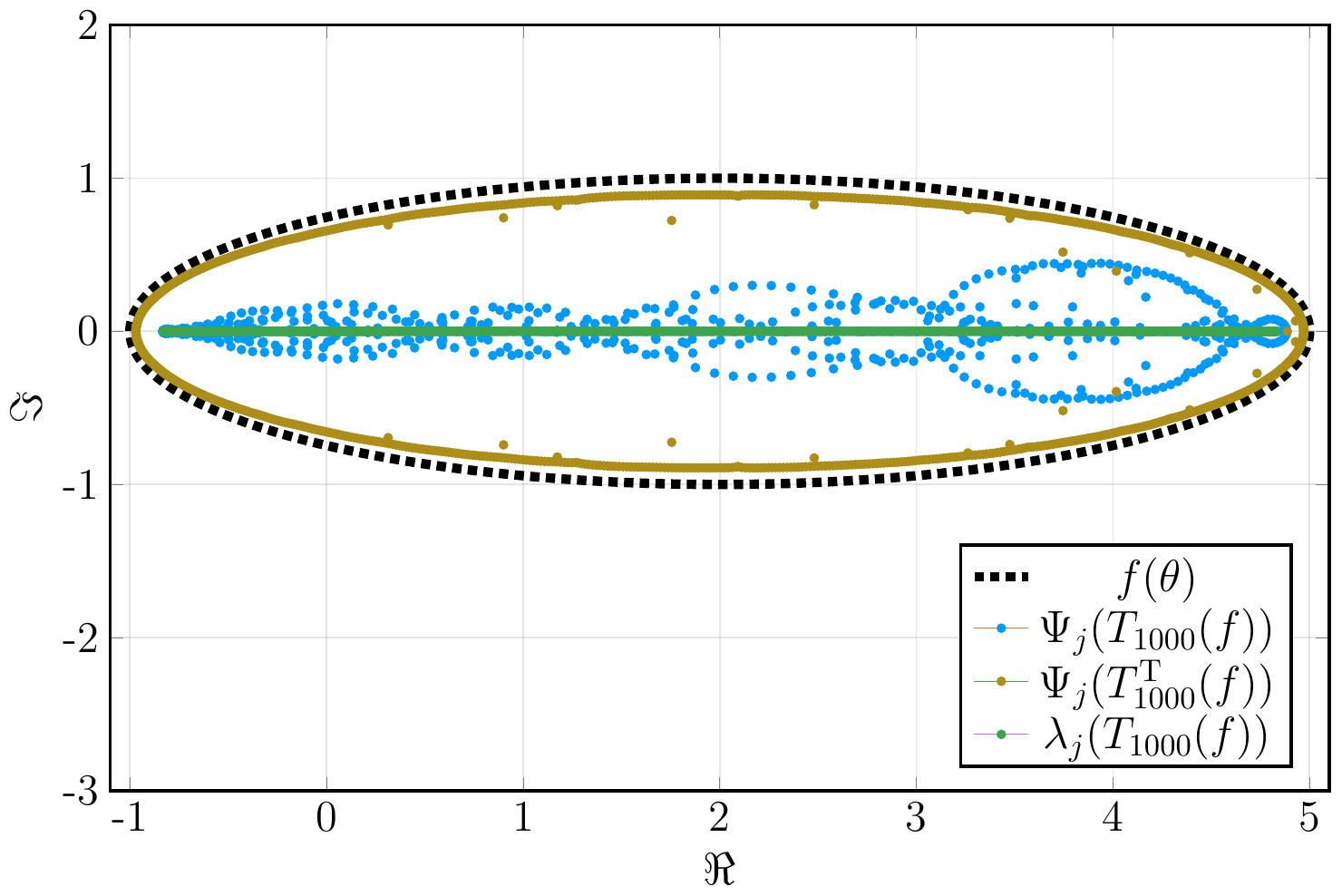}
\caption{[Example~\ref{exmp:1}: Symbol $f(\theta)=-e^{\mathbf{i}\theta}+2-2e^{-\mathbf{i}\theta}$]  Symbol $f(\theta)$ (dashed black line), the numerically computed spectra (using a standard double precision eigenvalue solver)  $\Psi_j(T_{1000}(f))$ (blue dots), $\Psi_j(T_{1000}^{\mathrm{T}}(f))$ (beige dots), and the analytical spectrum $\lambda_j(T_{1000}(f))=g(\theta_{j,1000})$ (green dots).}
\label{fig:exmp:1:pseudospectrum}
\end{figure}
In Figure~\ref{fig:exmp:1:pseudospectrum} we present the numerically computed spectra $\Psi_j(T_n(f))$ (blue dots) and $\Psi_j(T_n^{\mathrm{T}}(f))$ (beige dots), for $n=1000$, using a standard double precision eigenvalue solver. The analytical spectrum, defined by \eqref{eq:exmp:1:eigexacttridiag} and \eqref{eq:exmp:1:gexmptridiag} is also shown (green dots).
These numerically computed eigenvalues $\Psi_j(A_n)$ are related to the pseudospectrum, discussed for example in \cite{trefethen051,beam931,reichel921}. 
\end{exmp} 

\begin{exmp}
\label{exmp:2}
In this example we consider the symbol
\begin{align}
f(\theta)&=(2-2\cos(\theta))^2=6-8\cos(\theta)+2\cos(\theta)\nonumber
\end{align}
which generates a Toeplitz matrix $T_n(f)$ associated with the second order finite difference approximation of the bi-Laplacian, 
\begin{align}
T_n(f)=\left[
\begin{array}{rrrrrrrrrrr}
6&-4&1\vphantom{\ddots}\\
-4&6&-4&1\vphantom{\ddots}\\
1&-4&6&-4&1\vphantom{\ddots}\\
&\ddots&\ddots&\ddots&\ddots&\ddots\\
&&\ddots&\ddots&\ddots&\ddots&1\\
&&&\ddots&\ddots&\ddots&-4\\
&&&&1&-4&6\vphantom{\ddots}
\end{array}
\right].\nonumber
\end{align}
The matrices $T_n(f)$ are all Hermitian and sothey have a real spectrum. Moreover, we have $f(\theta)=g(\theta)$, and $\{T_n(f)\}_n\sim_{\sigma,\lambda}f$. In Figure~\ref{fig:exmp:2:spectrum} we represent the symbol $g=f$ and the eigenvalues of $T_n(f)$ for $n=5$.
The ``perfect'' sampling grid $\xi_{j,n}$ such that $\lambda_j(T_n(f))=g(\xi_{j,n})$ is not equispaced, but can in this case be obtained by either computing $\xi_j=2\sin^{-1}\left((\lambda_j(T_n(f)))^{1/4}/2\right)$ (since $f(\theta)=16\sin^4(\theta/2)$),
finding the roots in $(0,\pi)$ of $g(\theta)-\lambda_j(T_n(f))$, or using the expansion described in~\cite{ekstrom191} for large $n$.
\begin{figure}[!ht] 
\centering
\includegraphics[width=0.48\textwidth]{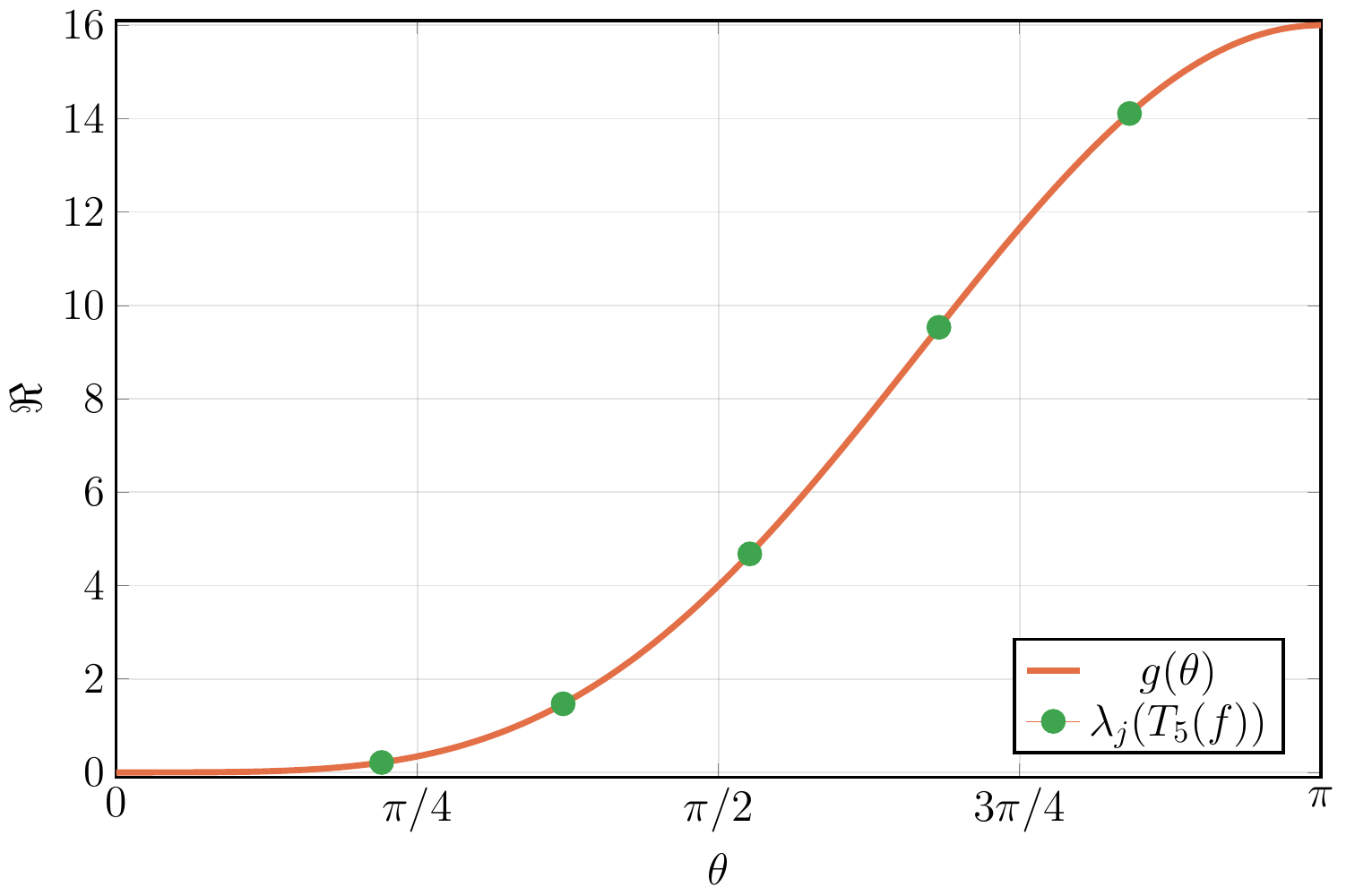}
\caption{[Example~\ref{exmp:2}: Symbol $f(\theta)=6-8\cos(\theta)+2\cos(2\theta)$] Symbol $g(\theta)=f(\theta)$ and $\lambda_j(T_5(f))=g(\xi_{j,5})$.}
\label{fig:exmp:2:spectrum}
\end{figure}
\end{exmp}

\begin{exmp}
\label{exmp:3}
In this example we consider the following symbol
\begin{align}
f(\theta)&=e^{\mathbf{i}\theta}-4+6e^{-\mathbf{i}\theta}-4e^{-2\mathbf{i}\theta}+e^{-3\mathbf{i}\theta}\nonumber\\
&=e^{-\mathbf{i}\theta}\left(6-8\cos(\theta)+2\cos(\theta)\right).\nonumber
\end{align}
The Toeplitz matrix $T_n(f)$ is a shifted version of the matrix considered in Example~\ref{exmp:2} (that is, the matrix associated with the second order finite difference approximation of the bi-Laplacian), and it is given by
\begin{align}
T_n(f)=\left[
\begin{array}{rrrrrrrrrrr}
-4&6&-4&1\vphantom{\ddots}\\
1&-4&6&-4&1\vphantom{\ddots}\\
&\ddots&\ddots&\ddots&\ddots&\ddots\\
&&\ddots&\ddots&\ddots&\ddots&1\\
&&&\ddots&\ddots&\ddots&-4\\
&&&&\ddots&\ddots&6\\
&&&&&1&-4\vphantom{\ddots}
\end{array}
\right].\nonumber
\end{align}
We note that
\begin{align}
f(\theta)&=e^{-\mathbf{i}\theta}\left(2-2\cos(\theta)\right)^2\nonumber\\
&=e^{-3\mathbf{i}\theta}\left(1-e^{\mathbf{i}\theta}\right)^4,\nonumber
\end{align} 
which is equivalent to (41) in \cite[Example 3.]{shapiro171} with $z=e^{\mathbf{i}\theta}$, $a=-1, r=3$, and $s=1$. Hence by (43) in the same article we have that
\begin{align}
g(\theta)=-\frac{\sin^4(\theta)}{\sin(\theta/4)\sin^3(3\theta/4)},
\label{eq:exmp:3:g}
\end{align}
and the matrix $T_n(g)$ would be full with $\lambda_j(T_n(f))\approx\lambda_j(T_n(g))\in(-\frac{(r+s)^{r+s}}{r^rs^s},0)=(-\frac{256}{27},0)$ for all $j$.

In the left panel of Figure~\ref{fig:exmp:3:symbols} we represent the functions $f$ (dashed black line), $g$ (red line) and the eigenvalues $\lambda_j(T_n(f))$ (green dots) for $n=5$. 
In the right panel of Figure~\ref{fig:exmp:3:symbols}  we show the function $g$ (red line) on the interval $[0,\pi]$ only (since it is even on $[-\pi,\pi]$) and the eigenvalues $\lambda_j(T_5(f))=g(\xi_{j,5})$ (green dots).
\begin{figure}[!ht] 
\centering
\includegraphics[width=0.47\textwidth,valign=t]{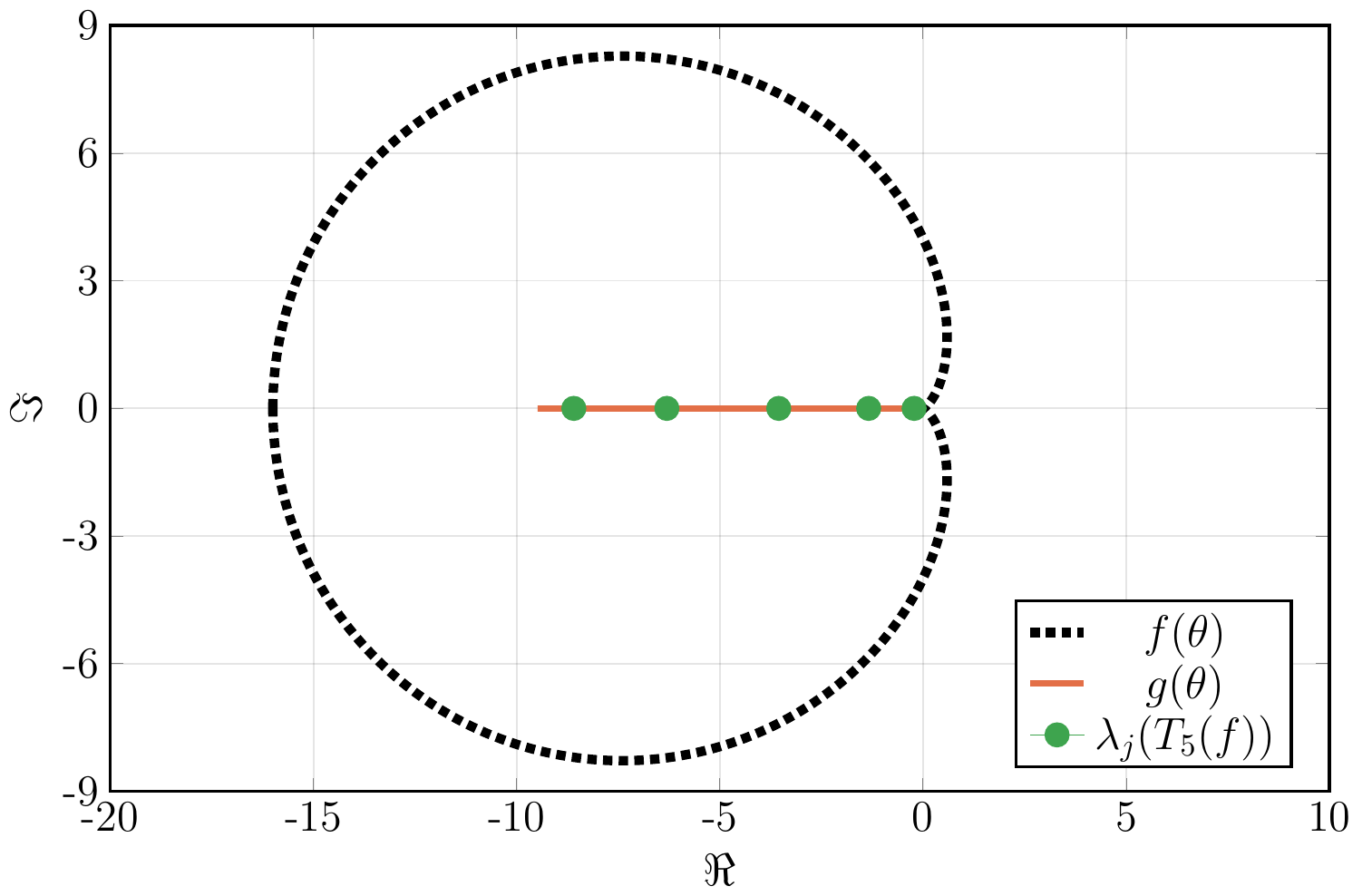}
\includegraphics[width=0.48\textwidth,valign=t]{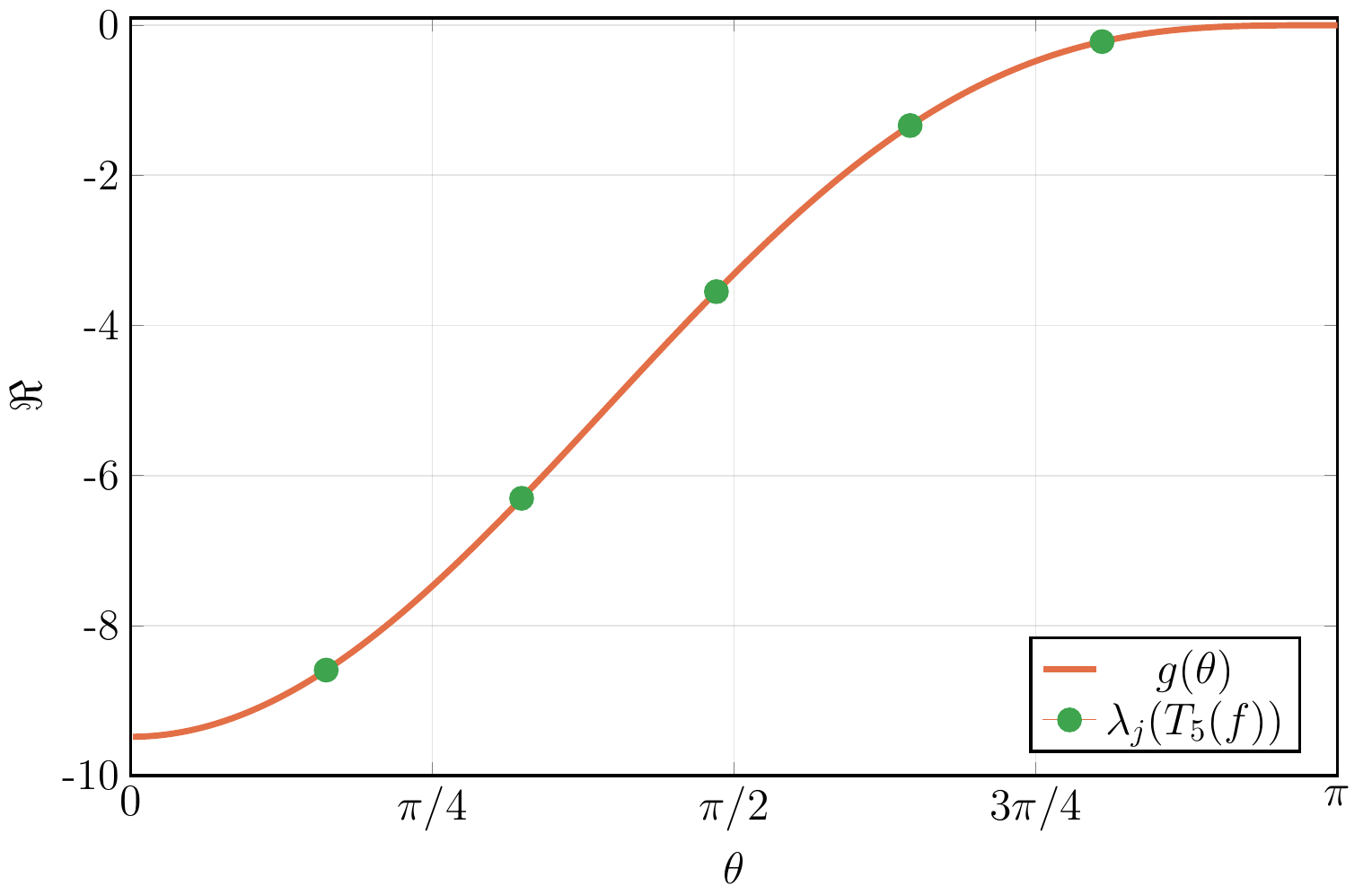}

\caption{[Example~\ref{exmp:3}: Symbol $f(\theta)=e^{-\mathbf{i}\theta}\left(6-8\cos(\theta)+2\cos(\theta)\right)$] Left: Representations of $f(\theta)$ (dashed black line),  $g(\theta)=-\sin^4(\theta)/(\sin(\theta/4)\sin^3(3\theta/4))$ (red line), and $\lambda_j(T_n(f))$ for $n=5$ (green dots). Right: Representation of $g$ and $\lambda_j(T_5(f))=g(\xi_{j,5})$.}
\label{fig:exmp:3:symbols}
\end{figure}
In Figure~\ref{fig:exmp:3:pseudospectrum} we present the numerically computed spectra $\Psi_j(T_n(f))$ (blue dots) and $\Psi_j(T_n^{\mathrm{T}}(f))$ (beige dots), for $n=1000$, using a standard double precision eigenvalue solver. The approximations of the true eigenvalues $\lambda_j(T_{1000}(f))=g{\xi_{j,1000}}$ (green dots) are also shown, computed with 128 bit precision. 
\begin{figure}[!ht] 
\centering
\includegraphics[width=0.48\textwidth]{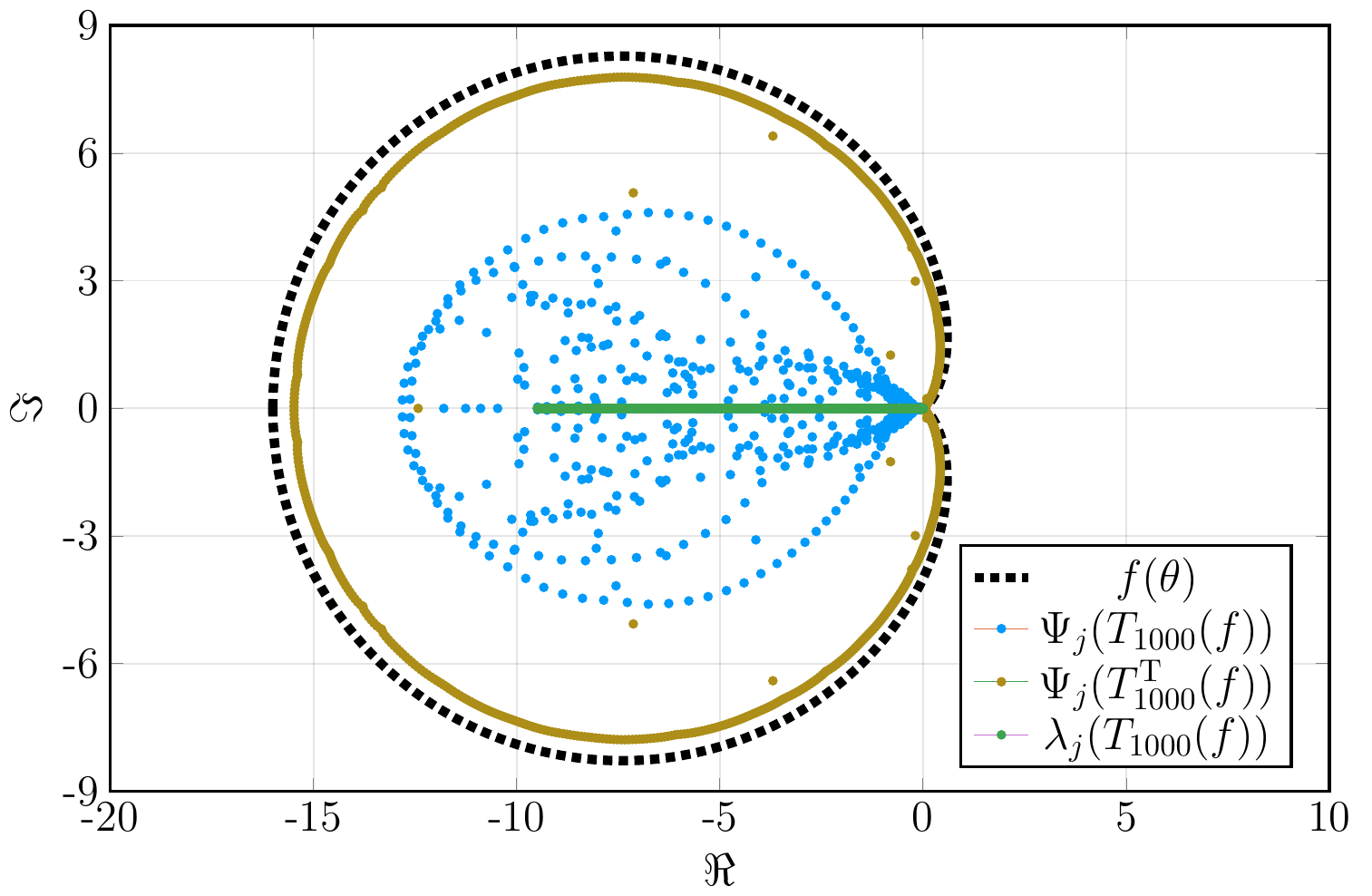}
\caption{[Example~\ref{exmp:3}: Symbol $f(\theta)=e^{-\mathbf{i}\theta}\left(6-8\cos(\theta)+2\cos(\theta)\right)$] Symbol $f(\theta)$ (dashed black line),
the numerically computed spectra (using a standard double precision eigenvalue solver) $\Psi_j(T_{1000}(f))$ (blue dots), $\Psi_j(T_{1000}^{\mathrm{T}}(f))$ (beige dots), and approximated spectrum $\lambda_j(T_{1000}(f))=g(\xi_{j,1000})$ (green dots).}
\label{fig:exmp:3:pseudospectrum}
\end{figure}

\end{exmp}

\begin{exmp}
\label{exmp:4}
In this example we consider a symbol $f$ where the explicit formula for $g$ is not known explicitly. Let
\begin{align}
f(\theta)=
 e^{3\mathbf{i}\theta}
-e^{2\mathbf{i}\theta}
+7e^{\mathbf{i}\theta}
+9e^{-1\mathbf{i}\theta}
-2e^{-2\mathbf{i}\theta}
+2e^{-3\mathbf{i}\theta}
-e^{-4\mathbf{i}\theta},\nonumber
\end{align}
which generates the matrix
\begin{align}
T_n(f)=\left[
\begin{array}{rrrrrrrrrrrrrrr}
0&9&-2&2&-1\vphantom{\ddots}\\
7&0&9&-2&2&-1\vphantom{\ddots}\\
-1&7&0&9&-2&2&-1\vphantom{\ddots}\\
1&-1&7&0&9&-2&2&-1\vphantom{\ddots}\\
&\ddots&\ddots&\ddots&\ddots&\ddots&\ddots&\ddots&\ddots\\
&&\ddots&\ddots&\ddots&\ddots&\ddots&\ddots&\ddots&-1\\
&&&\ddots&\ddots&\ddots&\ddots&\ddots&\ddots&2\\
&&&&\ddots&\ddots&\ddots&\ddots&\ddots&-2\\
&&&&&\ddots&\ddots&\ddots&\ddots&9\\
&&&&&&1&-1&\hphantom{-}7&0\vphantom{\ddots}\\
\end{array}
\right].\nonumber
\end{align}
From~\cite[Example 4.]{shapiro171} we have strong indications that approximately $\lambda_j(T_n(f))\in [-22.09,14.96]$ for all $j$.

In the left panel of Figure~\ref{fig:exmp:4:symbols} we represent the symbol $f$ (dashed black line) and the eigenvalues $\lambda_{1000}(T_n(f))$, computed with a 256 bit eigenvalue solver (dashed red line) since $g$ is not known. 
The eigenvalues $\lambda_j(T_n(f))$ (green dots) for $n=5$ are also shown. 
In the right panel of Figure~\ref{fig:exmp:4:symbols} show again the eigenvalues $\lambda_{1000}(T_n(f))$ arranged in non-decreasing order (dashed red line) since $g$ is not known on the interval $[0,\pi]$, since it is even on $[-\pi,\pi]$. Also the eigenvalues $\lambda_j(T_5(f))=g(\xi_{j,5})$ (green dots) are represented. The ``perfect'' grid $\xi_{j,n}$ is computed using data from Example~\ref{exmp:8}. Numerically we have $\lambda_{j}(T_{1000}(f))\in[-22.0912,14.9641]$ in agreement with~\cite{shapiro171}.

\begin{figure}[!ht] 
\centering
\includegraphics[width=0.472\textwidth,valign=t]{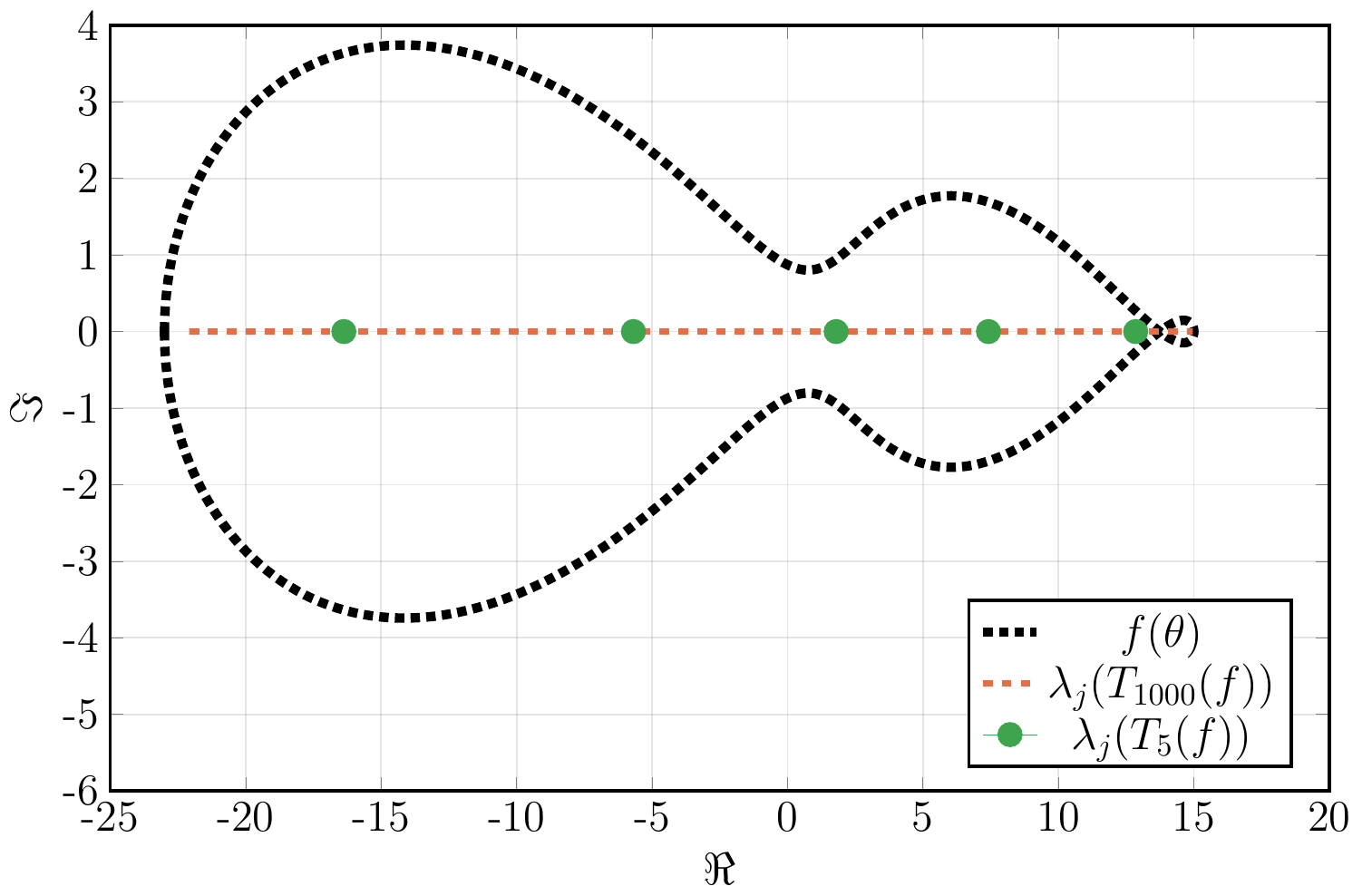}
\includegraphics[width=0.48\textwidth,valign=t]{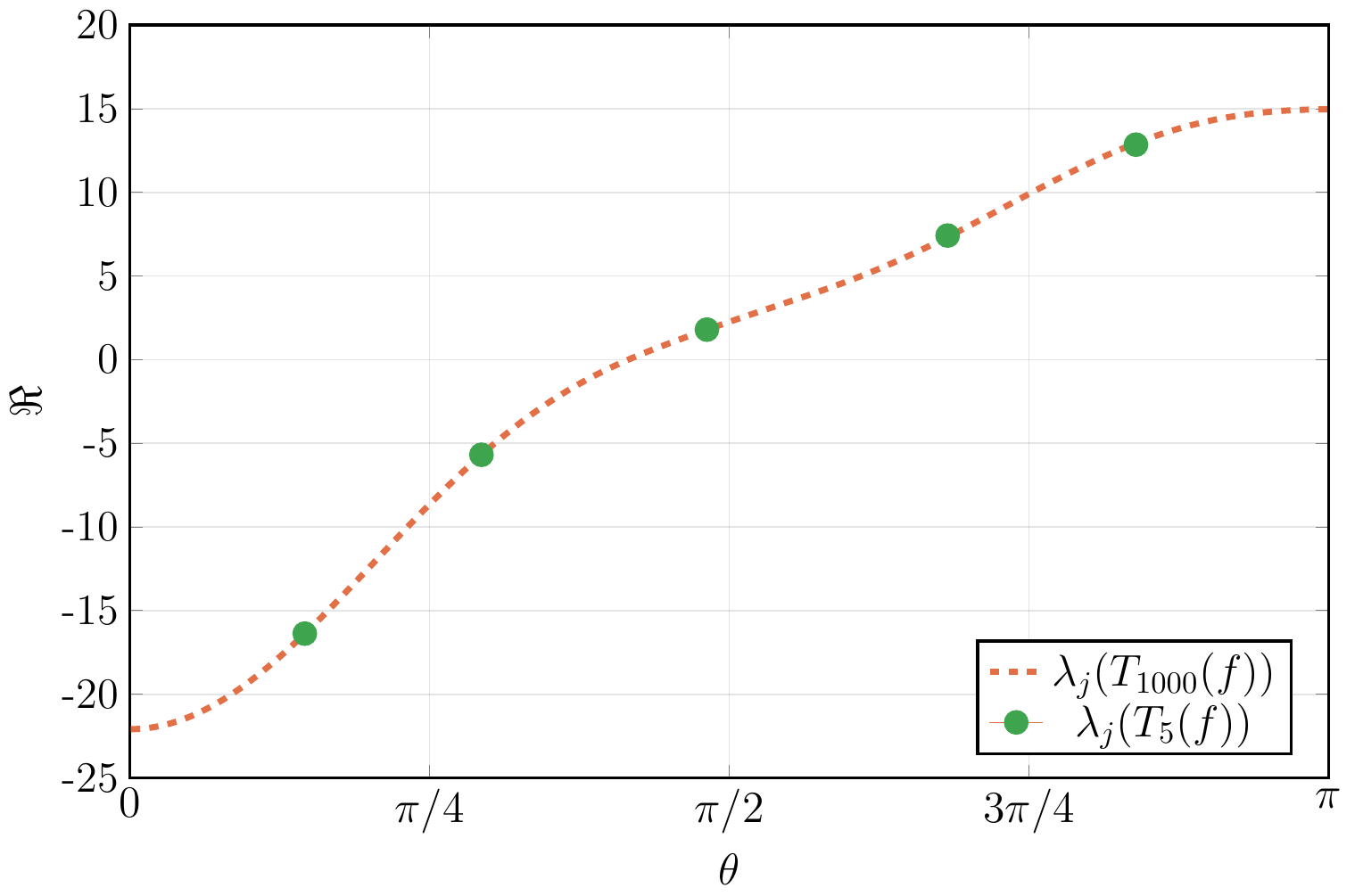}

\caption{[Example~\ref{exmp:4}: Symbol $f(\theta)=
 e^{3\mathbf{i}\theta}
-e^{2\mathbf{i}\theta}
+7e^{\mathbf{i}\theta}
+9e^{-1\mathbf{i}\theta}
-2e^{-2\mathbf{i}\theta}
+2e^{-3\mathbf{i}\theta}
-e^{-4\mathbf{i}\theta}$]  Left: Symbol $f(\theta)$ (dashed black line) and 
$\lambda_{j}(T_{1000}(f))$ (dashed red line) since $g(\theta)$ is unknown, and $\lambda_j(T_n(f))$ for $n=5$ (green dots).
Right: Eigenvalues $\lambda_{j}(T_{1000}(f))$ ordered in non-decreasing order and $\lambda_j(T_5(f))=g(\xi_{j,5})$.}
\label{fig:exmp:4:symbols}
\end{figure}
In Figure~\ref{fig:exmp:4:pseudospectrum} we present the numerically computed spectra $\Psi_j(T_n(f))$ (blue dots) and $\Psi_j(T_n^{\mathrm{T}}(f))$ (beige dots), for $n=1000$, using a standard double precision eigenvalue solver. The true eigenvalues $\lambda_j(T_{1000}(f))=g(\xi_{j,1000})$ (green dots) are approximated using a 256 bit precision computation.

\begin{figure}[!ht] 
\centering

\includegraphics[width=0.48\textwidth]{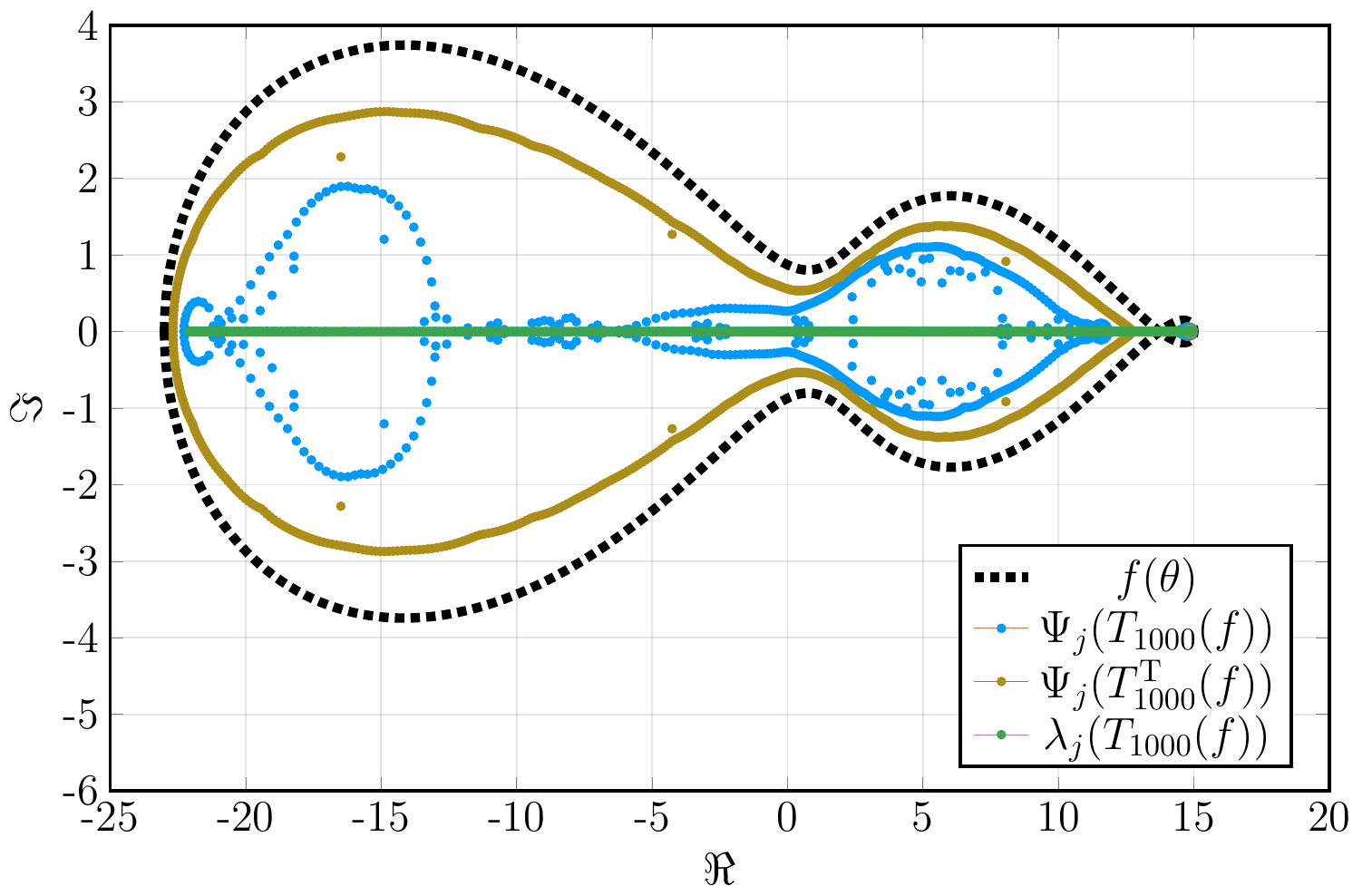}
\caption{[Example~\ref{exmp:4}: Symbol $f(\theta)=
 e^{3\mathbf{i}\theta}
-e^{2\mathbf{i}\theta}
+7e^{\mathbf{i}\theta}
+9e^{-1\mathbf{i}\theta}
-2e^{-2\mathbf{i}\theta}
+2e^{-3\mathbf{i}\theta}
-e^{-4\mathbf{i}\theta}$] Symbol $f(\theta)$ (dashed black line), the numerically computed spectra $\Psi_j(T_{1000}(f))$ (blue dots), $\Psi_j(T_{1000}^{\mathrm{T}}(f))$ (beige dots), and $\lambda_j(T_{1000}(f))=g(\xi_{j,n})$ (green dots).}
\label{fig:exmp:4:pseudospectrum}
\end{figure}

\end{exmp}

\section{Describing the real-valued eigenvalue distribution}
\label{sec:describing}

Assuming that $g$ is a real cosine trigonometric (RCTP) symbol associated with a symbol $f$ as in the working hypothesis, we introduce in Section~\ref{sec:describing:approximating} a new \textit{matrix-less} method to accurately compute the expansion functions $c_k, k=0,\ldots,\alpha$, where we recall that $c_0=g$.
Subsequently, in Section~\ref{sec:describing:constructing} we present procedures to obtain an approximation or even the analytical expression of $g$.

\subsection{Approximating the expansion functions $c_k$ in grid points $\theta_{j,n_0}$}
\label{sec:describing:approximating}

An asymptotic expansion of the eigenvalue errors $E_{j,n,0}\coloneqq E_{j,n}$, when sampling the symbol $f$ with the grid $\theta_{j,n}$ defined in the working hypothesis, under certain assumptions on $f$ implying that $g=f$, was discussed in a series of papers~\cite{bogoya151,bogoya171,bottcher101}; such expansion can be deduced from 
\begin{align}
\lambda_j(T_n(f))=f(\theta_{j,n})+\underbrace{\sum_{k=1}^\alpha c_k(\theta_{j,n})h^k+E_{j,n,\alpha}}_{E_{j,n}},\label{eq:describing:approximating:expl}
\end{align}
where $\theta_{j,n}$, $h$, and $E_{j,n,\alpha}$ are defined in the working hypothesis.

An algorithm was proposed in \cite{ekstrom171} to approximate the functions $c_k(\theta)$, which was subsequently extended to other types of Toeplitz-like matrices $A_n$ possessing an asymptotic expansion such as \eqref{eq:describing:approximating:expl}; see~\cite{ahmad171,ekstrom181,ekstrom185,ekstrom183}. We call this type of methods \textit{matrix-less}, since they do not need to construct the large matrix $A_n$ to approximate its eigenvalues; indeed, they approximate the functions $c_k(\theta)$ from $\alpha$ small matrices $A_{n_1},\ldots,A_{n_\alpha}$ and then they use this approximations to compute the approximate spectrum of~$A_n$ through the formula 
\begin{align}
\lambda_{j}(A_n)\approx \tilde{\lambda}_j(A_n)= f(\theta_{j,n})+\sum_{k=1}^\alpha \tilde{c}_k(\theta_{j,n})h^k.\label{eq:describing:approximating:eigvalapprox}
\end{align}
Assuming that the eigenvalues of $T_n(f)$ admit an asymptotic expansion in terms of an unknwn function $g$ instead of $f$, as in our working hypothesis, we can use a slight modification of Algorithm 1 in~\cite[Section~2.1]{ekstrom184} in order to find approximations of both $g$ and the eigenvalues of $T_n(f)$ through the following formula, analogous to \eqref{eq:describing:approximating:eigvalapprox}:
\begin{align}
\lambda_{j}(T_{n}(f))\approx \tilde{\lambda}_{j}(T_{n}(f))&=\sum_{k=0}^\alpha \tilde{c}_k(\theta_{j,n})h^k\nonumber\\
&=\tilde{g}(\theta_{j,n})+\sum_{k=1}^\alpha \tilde{c}_k(\theta_{j,n})h^k,
\label{eq:describing:approximating:expansion}
\end{align}
where the approximation $\tilde g(\theta)\coloneqq\tilde c_0(\theta)$ of $g(\theta)\coloneqq c_0(\theta)$ is obtained from $\alpha+1$ small matrices $T_{n_0}(f),\ldots,T_{n_\alpha}(f)$ as mentioned above. 

Here follows an implementation in \textsc{Julia} of the algorithm that computes the approximations $\tilde c_k(\theta)$ for $k=0,\ldots,\alpha$; the algorithm is written for clarity and not performance. All computations in this article are made with \textsc{Julia} 1.1.0~\cite{bezanson171}, using \texttt{Float64} and \texttt{BigFloat} data types, and the \textsc{GenericLinearAlgebra.jl} package~\cite{noack191}.
\begin{algo}
\label{algo:1}
Approximate expansion functions $c_k(\theta)$ for $k=0,\ldots,\alpha$ on the grid $\theta_{j,n_0}$.

\normalfont
{\footnotesize
\begin{lstlisting}[backgroundcolor = \color{jlbackground},
                   language = Julia,
                   xleftmargin = 0.1em,
                   framexleftmargin = 0.1em]
using LinearAlgebra, GenericLinearAlgebra 
setprecision(BigFloat,128)

# Example: computeC(100, 4, BigFloat[2, -1], BigFloat[2, -2])               
function computeC(n0        :: Integer,    # Number of grid points in grid theta_{j,n0}
                  alpha     :: Integer,    # Number of c_k to approximate, k=0,...,alpha
                  vc        :: Array{T,1}, # First column of T_n(f)
                  vr        :: Array{T,1}, # First row of T_n(f)
                  revorder  :: Bool=false  # Reverse ordering of eigenvalues of T_n(f)
                  ) where T 
  j0 = 1:n0
  E  = zeros(T,alpha+1,n0)
  hs = zeros(T,alpha+1)
  for kk = 1:alpha+1
    nk     = 2^(kk-1)*(n0+1)-1
    jk     = 2^(kk-1)*j0
    hs[kk] = convert(T,1)/(nk+1)
    Tnk    = Toeplitz(nk,vc,vr)
    eTnk   = eigvals(Tnk)
    if !isreal(eTnk)
      error("Spectrum not real. Decrease n0 or alpha, or use BigFloat with higher precision.")
    end
    eTnk    = sort(real.(eTnk),rev=revorder)
    E[kk,:] = eTnk[jk]
  end
  V = zeros(T,alpha+1,alpha+1)
  for ii = 1:alpha+1, jj = 1:alpha+1
    V[ii,jj] = hs[ii]^(jj-1)
  end
  return C=V\E # Output: Matrix C, size (alpha+1,n0), with approximations c_k(theta_{j,n0})
end

# Example: Toeplitz(100, Float64[2, -1], Float64[2, -2])
function Toeplitz(n  :: Integer,     # Order of Toeplitz matrix T_n(f)    
                  vc :: Array{T,1},  # First column of T_n(f)
                  vr :: Array{T,1}   # First row of T_n(f)
                  ) where T
  Tn = zeros(T,n,n)
  for ii = 1:length(vc)
    Tn = Tn + diagm(-ii+1 => vc[ii]*ones(T,n-ii+1))
  end
  for jj = 2:length(vr)
    Tn = Tn + diagm( jj-1 => vr[jj]*ones(T,n-jj+1))
  end
  return Tn # Output: Toeplitz matrix of order n, defined by vectors vc and vr
end
\end{lstlisting}
}
\end{algo}
Using the output $\tilde{c}_k(\theta_{j,n_0})$, we can employ the interpolation--extrapolation technique described in~\cite{ekstrom183} to efficiently compute very accurate approximations of $\tilde c_k(\theta)$ and, through \eqref{eq:describing:approximating:expansion}, the eigenvalues of $T_n(f)$ for an arbitrarily large order $n$. 
In the next section, we focus on the use of the approximations $\tilde{c}_0=\tilde g$ to describe $g$.

\subsection{Constructing a function $g$ from approximations $\tilde{g}(\theta_{j,n_0})=\tilde{c}_0(\theta_{j,n_0})$}
\label{sec:describing:constructing}
We here assume, for the sake of simplicity, that the sought function $g$ 
is real and even, so that it admits a cosine Fourier series of the form 
\begin{align}
g(\theta)=\hat{g}_0+2\sum_{k=1}^\infty \hat{g}_k\cos(k\theta),\qquad \hat{g}_k\in\mathbb{R}.
\label{eq:describing:constructing:g}
\end{align}
As we shall see, if $g$ is a real cosine trigonometric polynomial (RCTP), that is, a function of the form
\begin{align}
g(\theta)=\hat{g}_0+2\sum_{k=1}^m \hat{g}_k\cos(k\theta), \qquad \hat{g}_k\in\mathbb{R},
\label{eq:describing:constructing:gbanded}
\end{align}
then we will be able to recover the exact expression of $g$ (see Examples~\ref{exmp:5} and~\ref{exmp:6}); otherwise, we will get a truncated representation of the Fourier expansion of $g$ in \eqref{eq:describing:constructing:g} (see Examples~\ref{exmp:7} and~\ref{exmp:8}).
More specifically, what we do is the following: we consider the approximations $\tilde c_0(\theta_{j,n_0})$ provided by Algorithm~1 and we approximate the first $n_0$ Fourier coefficients $\hat g_0,\ldots,\hat g_{n_0}$ with the numbers $\tilde{\hat{g}}_0,\ldots,\tilde{\hat{g}}_{n_0}$ obtained by solving the linear system
\begin{align}
\tilde{\hat{g}}_0+2\sum_{k=1}^{n_0} \tilde{\hat{g}}_k\cos(k\theta_{j,n_0})=\tilde c_0(\theta_{j,n_0}),\qquad j=1,\ldots,n_0.\label{eq:describing:constructing:tildehatgsystem}
\end{align}

\begin{algo}
\label{algo:2}
Compute approximations $\tilde{\hat{g}}_k$ of the Fourier coefficients $\hat{g}_k$ of $g(\theta)$.
\normalfont
{\footnotesize
\begin{lstlisting}[backgroundcolor = \color{jlbackground},
                   language = Julia,
                   xleftmargin = 0.1em,
                   framexleftmargin = 0.1em]
# Example: computeghattilde(C[1,:])
function computeghattilde(c0 :: Array{T,1}) where T # Array of approximations c_0(theta_{j,n0})
  n0 = length(c0)
  t = LinRange(convert(T,pi)/(n0+1),n0*convert(T,pi)/(n0+1),n0)
  G = zeros(T,n0,n0)
  G[:,1] = ones(T,n0)
  for jj = 2:n0
    G[:,jj] = 2*cos.((jj-1)*t)
  end
  return ghattilde = G\c0 # Output: Coefficients ghattilde in (13)
end
\end{lstlisting}
}
\end{algo}

\section{Numerical examples}
\label{sec:numerical}
We now employ the proposed Algorithms~\ref{algo:1} and \ref{algo:2} on a the symbols $f$ discussed in Examples~\ref{exmp:1}--\ref{exmp:4} to highlight the applicability of the approach, in the respective Examples~\ref{exmp:5}--\ref{exmp:8}.

\begin{itemize}
\item Example~\ref{exmp:5}: Only $\tilde{c}_0$ is non-zero, since $\theta_{j,n}$ gives exact eigenvalues, and the function $g$ is constructed.
\item Example~\ref{exmp:6}: Symbol $f=c_0$, and $c_k,k=1,\ldots,4$, are recovered accurately, and the function $g=f$ is constructed. 
\item Example~\ref{exmp:7}: Symbol $g=c_0$, and $c_k,k=1,\ldots,4$, are recovered accurately, and a truncated RCTP representation of of $g$ is constructed.
\item  Example~\ref{exmp:8}: Symbol $g=c_0$, and $c_k,k=1,\ldots,4$, are constructed, and a truncated RCTP representation of of $g$ is constructed.
\end{itemize}
\begin{exmp}
\label{exmp:5}
We return to the non-symmetric symbol $f(\theta)=-e^{\mathbf{i}\theta}+2-2e^{-\mathbf{i}\theta}$ of Example~\ref{exmp:1}, and first use the proposed Algorithm~\ref{algo:1}.
Note that care has to be taken when using for example standard \textsc{Matlab} \texttt{eig} command, since already for $n=160$ the returned eigenvalues are complex-valued (and wrong). In such a circumstance a choice of $n_0$ and $\alpha$ needs to be such that $n_\alpha=2^\alpha(n_0+1)-1<160$. However, we here also use an arbitrary precision solver, \textsc{GenericLinearAlgebra.jl} in \textsc{Julia}, so we can increase precision such that theoretically any combination of $n_0$ and $\alpha$ can be chosen. The performance however decreases fast as we increase the computational precision, showing the need for the current proposed algorithms.

In Figure~\ref{fig:exmp:5:expansion} we present the computation of $\tilde{c}_k(\theta_{j,n_0})$, for $n_0=31$, and different precision and $\alpha$.
In the left panel we show the approximated expansion functions $\tilde{c}_k$, $k=0,\ldots,\alpha$, where $\alpha=4$, and 128 bit precision computation.
As is seen, the only non-zero $\tilde{c}_k$ is $\tilde{c}_0$, which is expected since the exact eigenvalues are given by $g(\theta_{j,n_0})=c_0(\theta_{j,n})$.
In the right panel of Figure~\ref{fig:exmp:5:expansion} we show the absolute error in the approximation of $c_0(\theta_{j,n_0})$ for double precision computation, with $\alpha=2$, and 128 bit computation, with $\alpha=4$.

\begin{figure}[!ht]
\centering
\includegraphics[width=0.472\textwidth,valign=t]{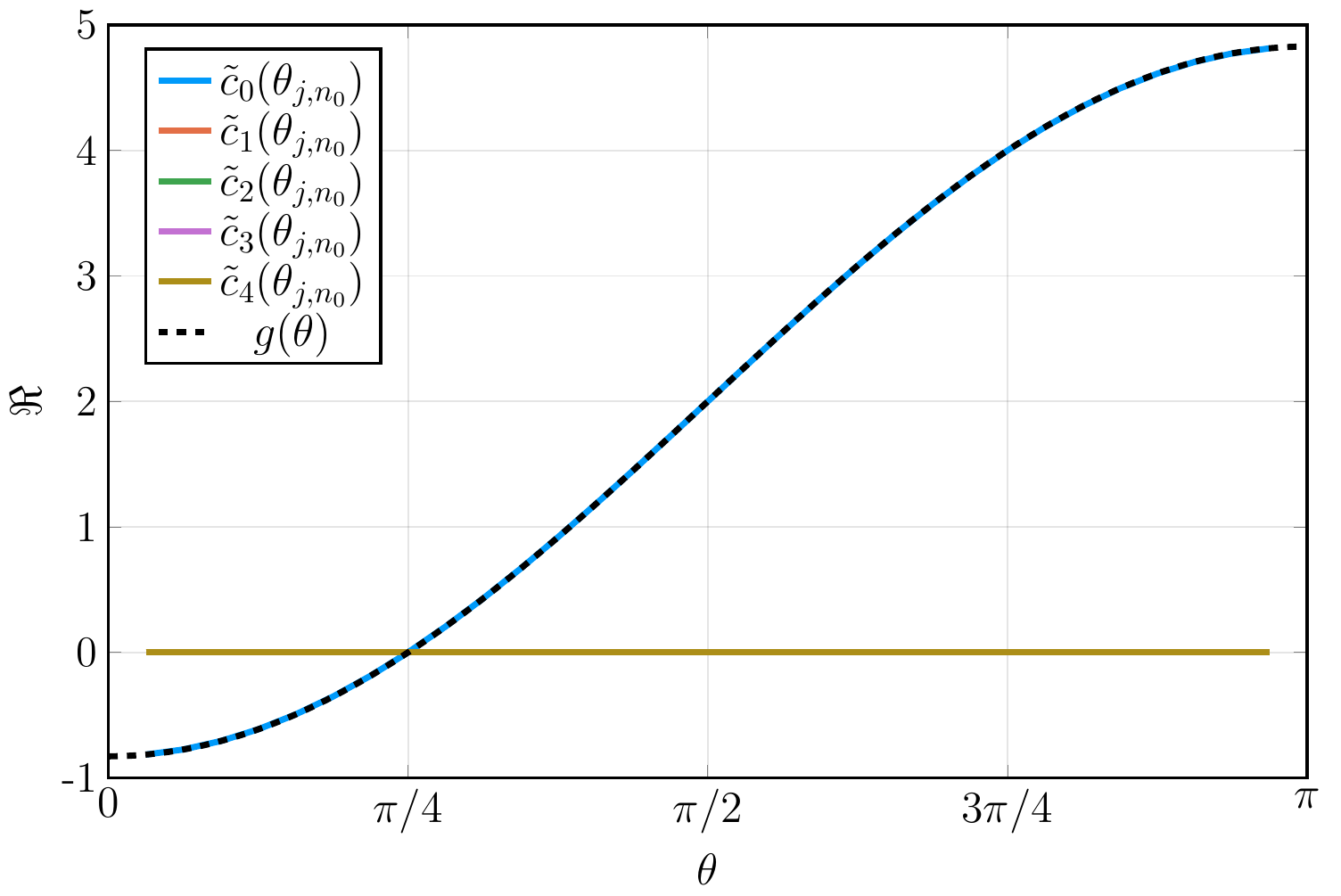}
\includegraphics[width=0.48\textwidth,valign=t]{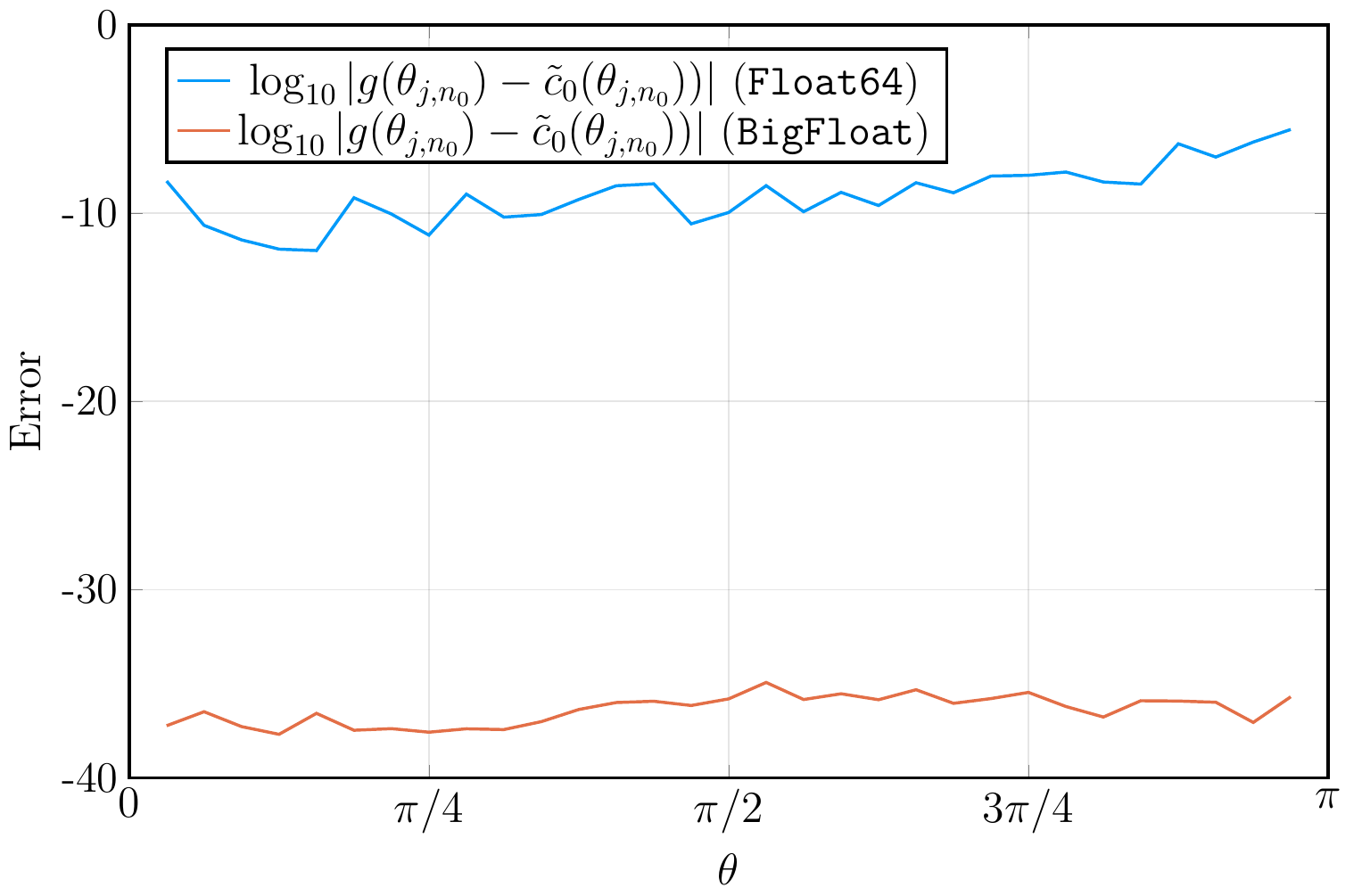}
\caption{[Example~\ref{exmp:5}: Symbol $f(\theta)=-e^{\mathbf{i}\theta}+2-2e^{-\mathbf{i}\theta}$] 
The computed $\tilde{c}_k(\theta_{j,n_0})$, $k=0,\ldots\alpha$, $n_0=31$ using Algorithm~\ref{algo:1}.  
Left: Computation using 128 bit precision with $\alpha=4$. Only $\tilde{c}_0(\theta_{j,n_0})=\tilde{g}(\theta_{j,n_0})$ are non-zero. Right: The absolute errors of $\tilde{c}_0(\theta_{j,n})$ compared with $g(\theta_{j,n})$ for double precision computation, with $\alpha=2$, and 128 bit precision, with $\alpha=4$.}
\label{fig:exmp:5:expansion}
\end{figure}

Now we employ Algorithm~\ref{algo:2} to compute the Fourier coefficients of $g$. For illustrative purposes we only show a subset of the system \eqref{eq:describing:constructing:tildehatgsystem}.
If we choose $n_0=2^\beta-1$, we have three grid points $\theta_{j_\beta,n_0}$ equal to $\pi/4, \pi/2$ and $3\pi/4$, corresponding to indices $j_\beta=2^{\beta-2}\{1,2,3\}$. We then have 
\begin{align}
\left[
\begin{array}{ccc}
1&\sqrt{2}&0\\
1&0&-1\\
1&-\sqrt{2}&0
\end{array}
\right]
\left[\begin{array}{c}
\tilde{\hat{g}}_0\\
\tilde{\hat{g}}_1\\
\tilde{\hat{g}}_2
\end{array}
\right]\approx
\left[\begin{array}{c}
\tilde{c}_0(\pi/4)\\
\tilde{c}_0(\pi/2)\\
\tilde{c}_0(3\pi/4)
\end{array}
\right].\nonumber
\end{align}
We construct the following system with $\tilde{c}_0$ computed with $n_0=31$ (that is, $\beta=5$ and $\alpha=2$ above) using double precision, and subsequently we compute $[\tilde{\hat{g}}_0, \tilde{\hat{g}}_1, \tilde{\hat{g}}_2]^{\mathrm{T}}$,
\begin{align}
\left[
\begin{array}{ccc}
1&\sqrt{2}&0\\
1&0&-1\\
1&-\sqrt{2}&0
\end{array}
\right]
\left[\begin{array}{c}
\tilde{\hat{g}}_0\\
\tilde{\hat{g}}_1\\
\tilde{\hat{g}}_2
\end{array}
\right]\approx
\left[\begin{array}{r}
 0.000000000006924\\
 1.999999999889880\\
 3.999999989588136
\end{array}
\right],\quad
\left[\begin{array}{c}
\tilde{\hat{g}}_0\\
\tilde{\hat{g}}_1\\
\tilde{\hat{g}}_2
\end{array}
\right]\approx
\left[\begin{array}{r}
  1.999999994797530\\
 -1.414213558689498\\
 -0.000000005092350
\end{array}
\right].\nonumber
\end{align}
We conclude from this computation that $g(\theta)=\hat{g}_0+2\hat{g}_1\cos(\theta)+2\hat{g}_2\cos(2\theta)=2-2\sqrt{2}\cos(\theta)$, which is the already known analytical expression; see \eqref{eq:exmp:1:eigexacttridiag} and \eqref{eq:exmp:1:gexmptridiag}. Note also that, in this simple example, the vector containing $\tilde{c}_0(\theta_{j_\beta,n_0})$ can be assumed to be equal to $[0,2,4]^{\mathrm{T}}$, which would yield the exact solution (to machine precision).
Using the full system \eqref{eq:describing:constructing:tildehatgsystem} in Algorithm~\ref{algo:2} yields the same result. If we now would approximate the monotonically non-increasing $g$ (instead of the non-decreasing) in Algorithm~\ref{algo:1} the vector containing $\tilde{c}(\theta_{j_\beta,n_0})$ would be $[4,2,0]^{\mathrm{T}}$ and would yield the symbol $g(\theta)=2+2\sqrt{2}\cos(\theta)$. Obviously, the eigenvalues of $T_n(g)$ are the same for both versions of $g$.
\end{exmp}

\begin{exmp}
\label{exmp:6}
We here return to the symmetric symbol $f(\theta)=(2-2\cos(\theta))^2=6-8\cos(\theta)+2\cos(2\theta)$, as in Example~\ref{exmp:2}.
Since we know
$\{T_n(f)\}_n\sim_{\sigma,\lambda}f=g$,
employing Algorithm~\ref{algo:1} will return as $\tilde{c}_0$ an approximation of $g$, and as $\tilde{c}_k$, $k>0$ the expansion functions previously obtained and studied in \cite{barrera181,ekstrom184}.

In Figure~\ref{fig:exmp:6:expansion} we show in the left panel the approximated expansion functions $\tilde{c}_k(\theta_{j,n_0})$ for $k=0,\ldots,\alpha$, computed using $n_0=100$, $\alpha=4$. Computations are made with double precision.
In the right panel of Figure~\ref{fig:exmp:6:expansion} we show the absolute error of the approximation of $g$, that is,  $\log_{10}|g(\theta_{j,n_0})-\tilde{c}_0(\theta_{j,n_0})|$.

The erratic behavior of $\tilde{c}_4(\theta)$ close to $\theta=0$ in the left panel, and the increased error close to to $\theta=0$ in the right panel are due to the fact that the symbol $f$ violates the so-called simple-loop conditions, discussed in  \cite{barrera181,ekstrom184}.
\begin{figure}[!ht]
\centering
\includegraphics[width=0.472\textwidth,valign=t]{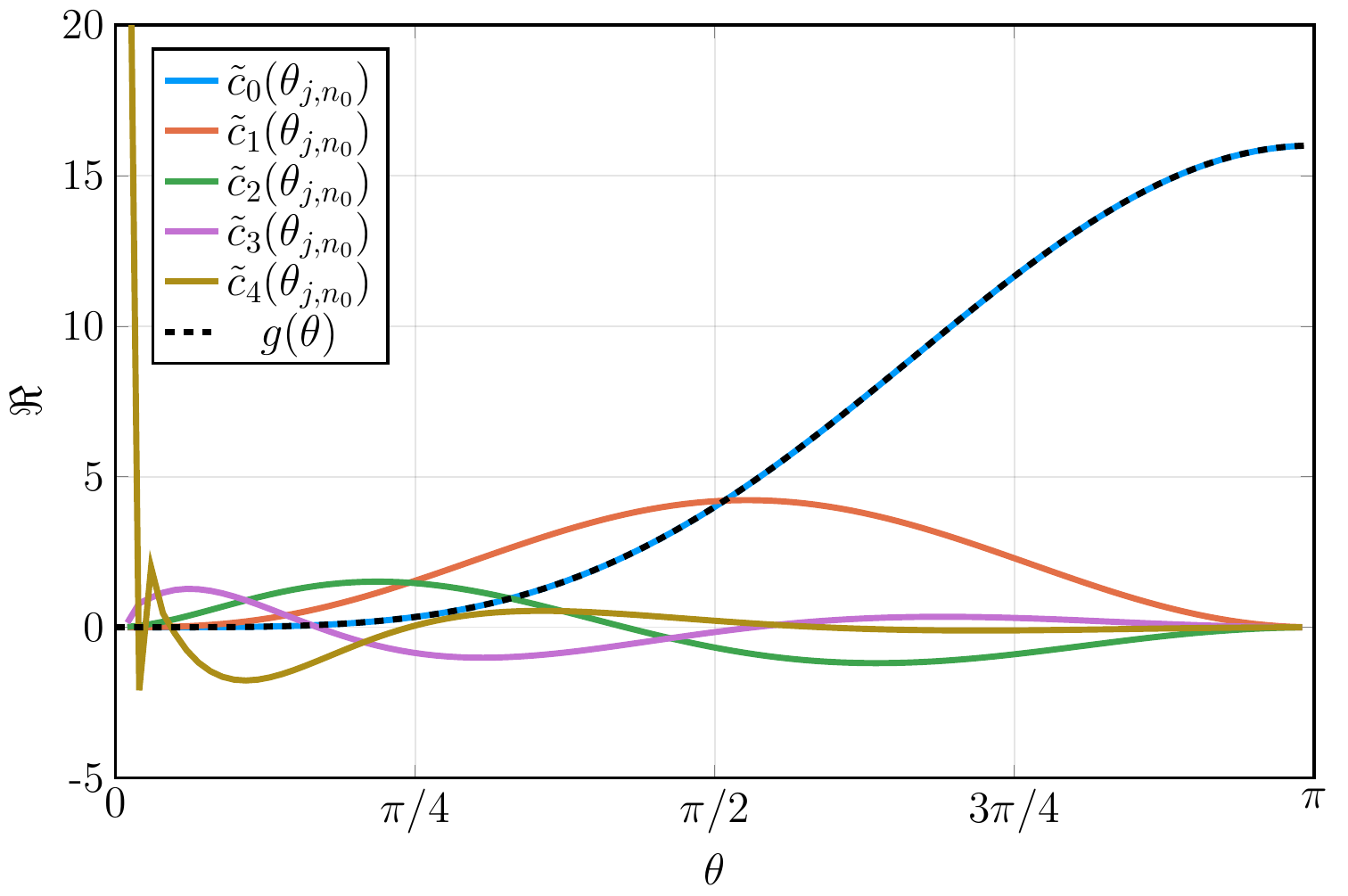}
\includegraphics[width=0.48\textwidth,valign=t]{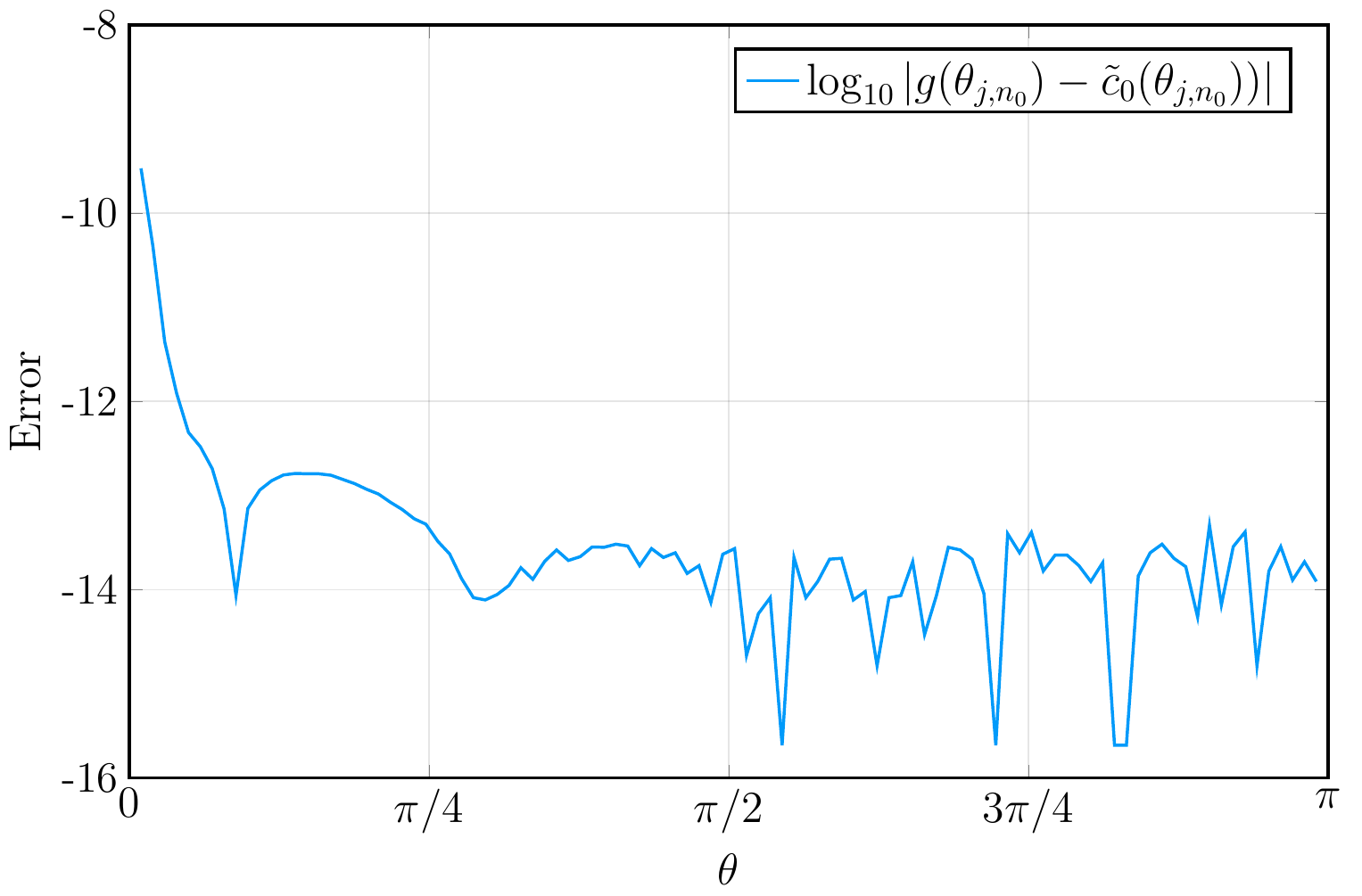}

\caption{[Example~\ref{exmp:6}: Symbol $f(\theta)=6-8\cos(\theta)+2\cos(2\theta)$]
 Left: The approximated expansion functions $\tilde{c}_k(\theta_{j,n_0}), k=0,\ldots,\alpha$ for $n_0=100$ and $\alpha=4$. Right: The absolute error $\log_{10}|g(\theta_{j,n_1})-\tilde{c}_0(\theta_{j,n_1})|$.}
\label{fig:exmp:6:expansion}
\end{figure}

Using Algorithm~\ref{algo:2} we compute approximations of the Fourier coefficients of the mononically increasing $g$ to be $\tilde{\hat{g}}_0=6,\tilde{\hat{g}}_1=-4$, $\tilde{\hat{g}}_2=1$, and $\tilde{\hat{g}}_k=0$ for $k>2$. We thus recover the true symbol $g(\theta)=f(\theta)=\hat{g}_0+2\hat{g}_1\cos(\theta)+2\hat{g}_2\cos(2\theta).$ If Algorithm~\ref{algo:1} was used to compute $\tilde{c}_k$ for the monotonically decreasing $g$ instead, the computed Fourier coefficients would be  $\tilde{\hat{g}}_0=6,\tilde{\hat{g}}_1=4$, and $\tilde{\hat{g}}_2=1$. In fact, for $g(\theta)=6\pm 8\cos(\theta)+2\cos(2\theta)$ we have the same eigenvalues for $T_n(f)$ and $T_n(g)$.
\end{exmp}

\begin{exmp}
\label{exmp:7}
In this example we continue the investigation of $f(\theta)=e^{-\mathbf{i}\theta}(6-8\cos(\theta)+2\cos(2\theta))$ from Example~\ref{exmp:3}. In Figure~\ref{fig:exmp:7:expansion} we show in the left panel the approximated expansion functions $\tilde{c}_k(\theta_{j,n_0})$ for $n_0=100$ and $\alpha=4$. Computations are made with 256 bit precision and the approximation $\tilde{c}_0(\theta_{j,n_0})$ overlaps well with $g$, defined in \eqref{eq:exmp:3:g}. Note the erratic behavior of $\tilde{c}_4$ close to $\theta=\pi$.
In the right panel of Figure~\ref{fig:exmp:7:expansion} we show the absolute values of the first one houndred approximated Fourier coefficients $\tilde{\hat{g}}_k$, given by Algorithm~\ref{algo:2}. 
In Table~\ref{tbl:exmp:7} we present the first ten true Fourier coefficients, $\hat{g}_k$, computed with $g$ defined in \eqref{eq:exmp:3:g} and \eqref{eq:introduction:fourier}, and the approximations $\tilde{\hat{g}}_k$ from Algorithm~\ref{algo:2}. Since $g$ is not an RCTP we can not recover the original simple expression of the symbol \eqref{eq:exmp:3:g}, but we can anyway obtain an approximated expression of $g$ through our Algorithm~\ref{algo:2}. 
\begin{figure}[!ht]
\centering
\includegraphics[width=0.475\textwidth,valign=t]{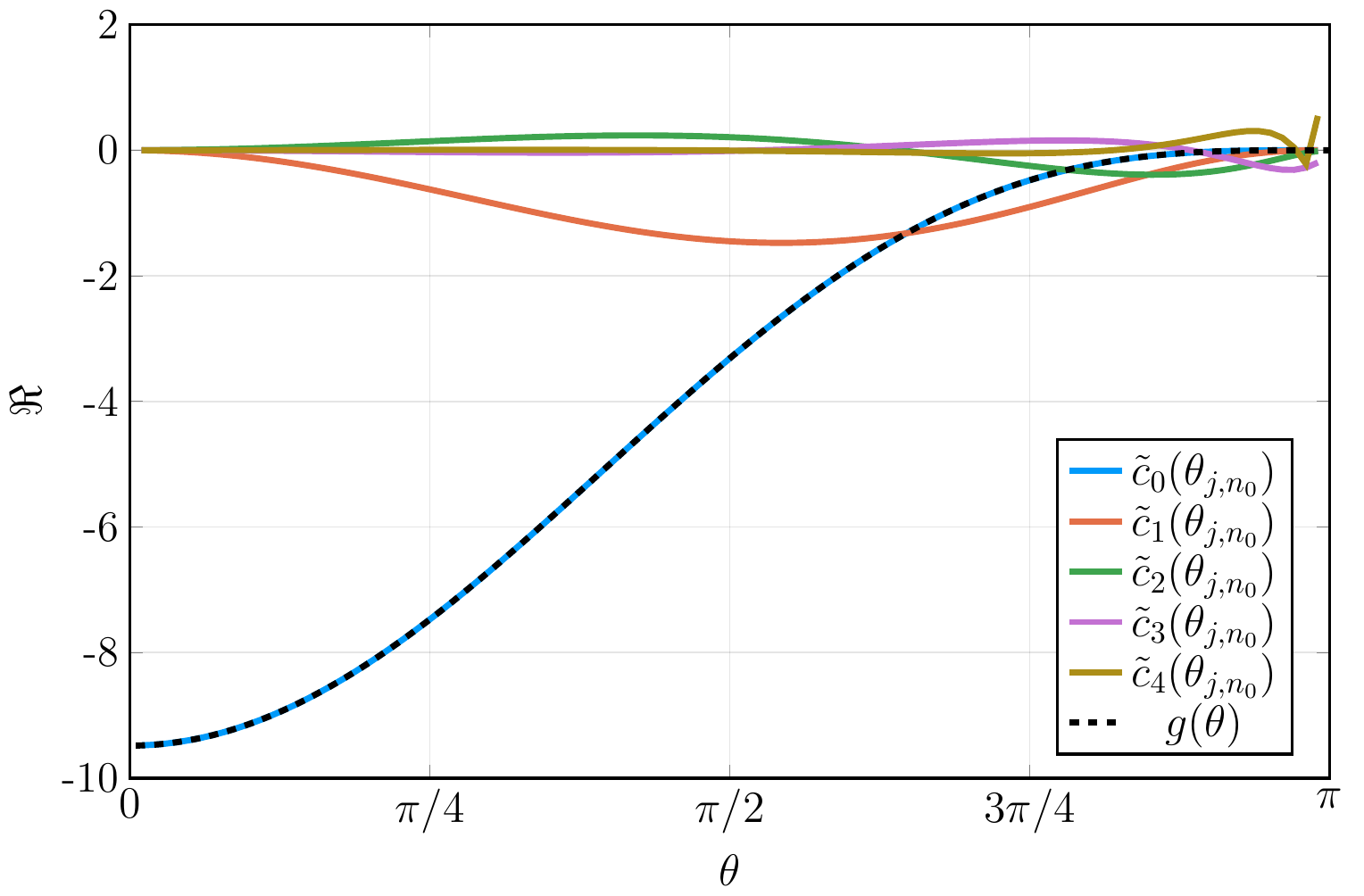}
\includegraphics[width=0.48\textwidth,valign=t]{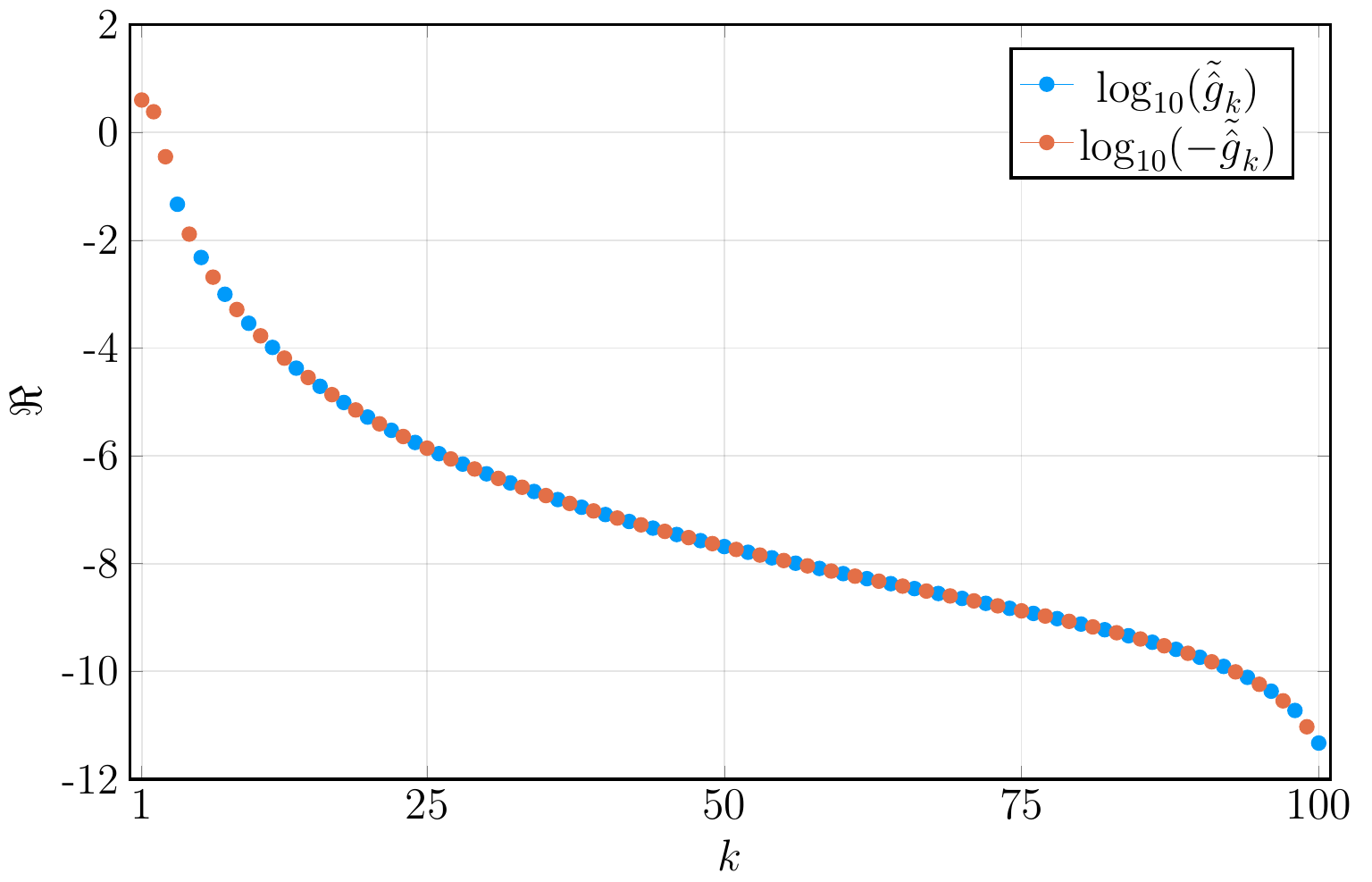}

\caption{[Example~\ref{exmp:7}: Symbol $f(\theta)=e^{-\mathbf{i}\theta}(6-8\cos(\theta)+2\cos(2\theta))$] Left: The approximated expansion functions $\tilde{c}_k(\theta_{j,n_0}), k=0,\ldots,\alpha$ for $n_0=100$ and $\alpha=4$.
The approximation $\tilde{c}_0(\theta_{j,n_0})$ overlaps well with $g(\theta)=-\sin^4(\theta)/(\sin(\theta/4)\sin^3(3\theta/4))$. Note the erratic behaviour of $\tilde{c}_4$ close to $\theta=\pi$. Right: The absolute value of the approximated first one houndred Fourier coefficients, $\log_{10}|\tilde{\hat{g}}_k|$. }
\label{fig:exmp:7:expansion}
\end{figure}

\begin{table}[!ht]
\centering
\caption{[Example~\ref{exmp:7}: Symbol $f(\theta)=e^{-\mathbf{i}\theta}(6-8\cos(\theta)+2\cos(2\theta))$] First ten true ($\hat{g}_k$) and computed ($\tilde{\hat{g}}_k$) Fourier coefficients of $g$. Approximations computed using $n_0=100$ and $\alpha=4$, and 256 bit precision.}
\label{tbl:exmp:7}
\begin{tabular}{rrrr}
\toprule
$k$ &$\hat{g}_k$&$\tilde{\hat{g}}_k$\\
\midrule
0&  $-4.000000000000000$&  $-3.999999999436239$\\
1&  $-2.423215805461417$& $-2.423215806024005$\\
2&  $-0.354481702999765$& $-0.354481702436023$\\
3&    $0.046583829909932$& $0.046583829347381$\\
4&   $-0.013008232443064$&$-0.013008231879376$\\
5&   $ 0.004790313798591$& $0.004790313236114$\\
6&  $-0.002068441503570$&  $-0.002068440939976$\\
7&   $ 0.000995276400689$& $0.000995275838326$\\
8&  $-0.000518988396995$& $-0.000518987833535$\\
9&   $0.000288215823752$&  $0.000288215261541$\\
\bottomrule
\end{tabular}
\end{table}
\end{exmp}

\begin{exmp}
\label{exmp:8}
Finally, we return to the non-symmetric symbol discussed in Example~\ref{exmp:4}, that is, $f(\theta)=
 e^{3\mathbf{i}\theta}
-e^{2\mathbf{i}\theta}
+7e^{\mathbf{i}\theta}
+9e^{-1\mathbf{i}\theta}
-2e^{-2\mathbf{i}\theta}
+2e^{-3\mathbf{i}\theta}
-e^{-4\mathbf{i}\theta}$.
Again, we employ Algorithms~\ref{algo:1} and \ref{algo:2} to study the symbols $f$ and $g$. In the left panel of Figure~\ref{fig:exmp:8:expansion} we present the approximated expansion functions in the working hypothesis, for $n_0=100$ and $\alpha=4$. Computations are made with 512 bit precision. The blue line, $\tilde{c}_0(\theta_{j,n_0})$ corresponds to the approximation of the unknown symbol $g$. Recall the curve of $\lambda_{j}(T_{1000}(f))$ in the right panel of Figure~\ref{fig:exmp:4:symbols}, which in principal matches the current $\tilde{c}_0$.
Note how all $\tilde{c}_k$, for $k>0$, are zero in apparently the same point $\theta_0\in[\frac{55\pi}{101},\frac{56\pi}{101}]$.
In the right panel of Figure~\ref{fig:exmp:8:expansion} we see the first one houndred approximated Fourier coefficients of $g$, by using Algorithm~\ref{algo:2}. In Table~\ref{tbl:exmp:8} is presented the first ten approximated Fourier coefficients, $\tilde{\hat{g}}_k$. 

\begin{figure}[!ht]
\centering

\includegraphics[width=0.475\textwidth,valign=t]{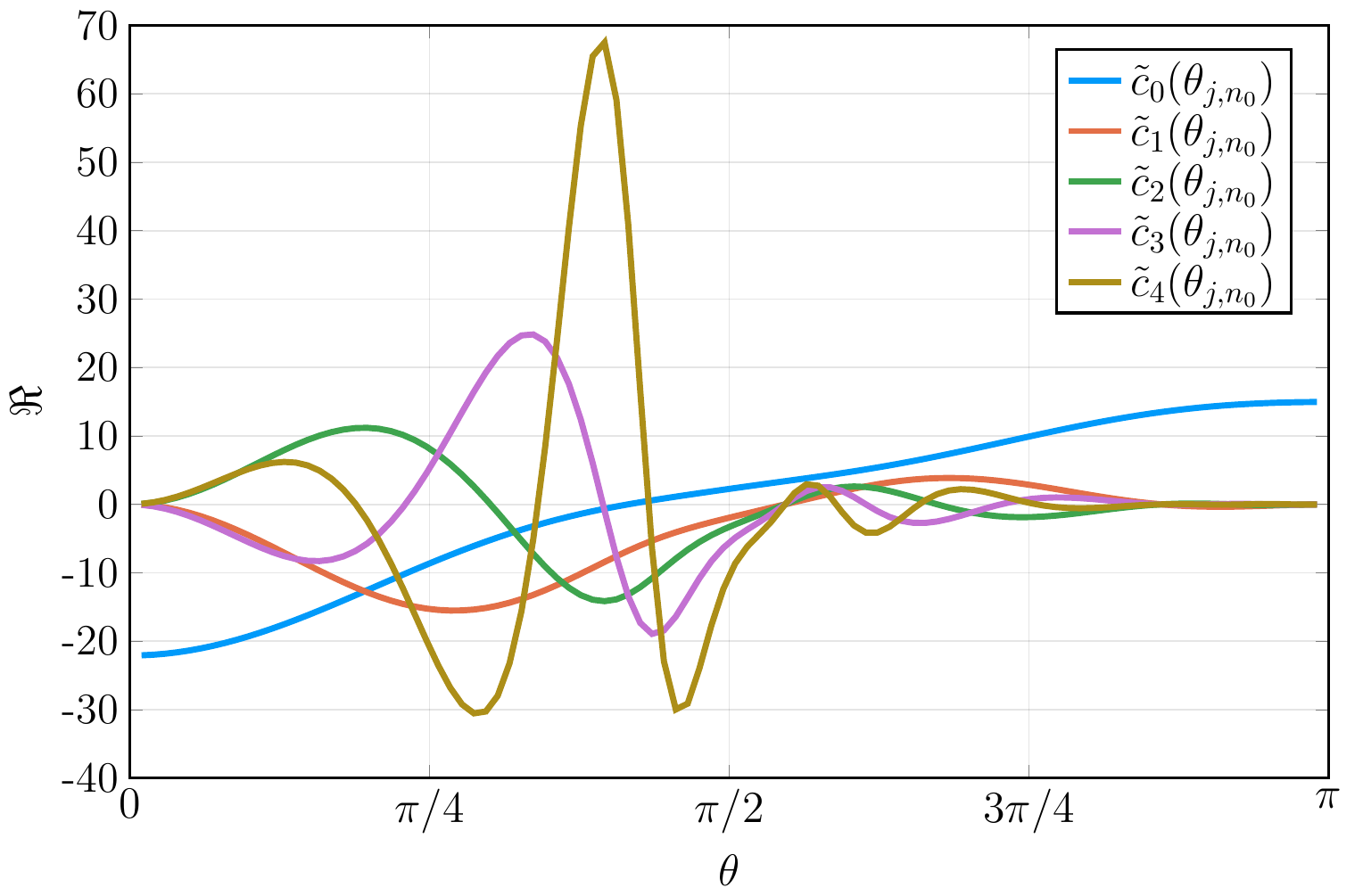}
\includegraphics[width=0.48\textwidth,valign=t]{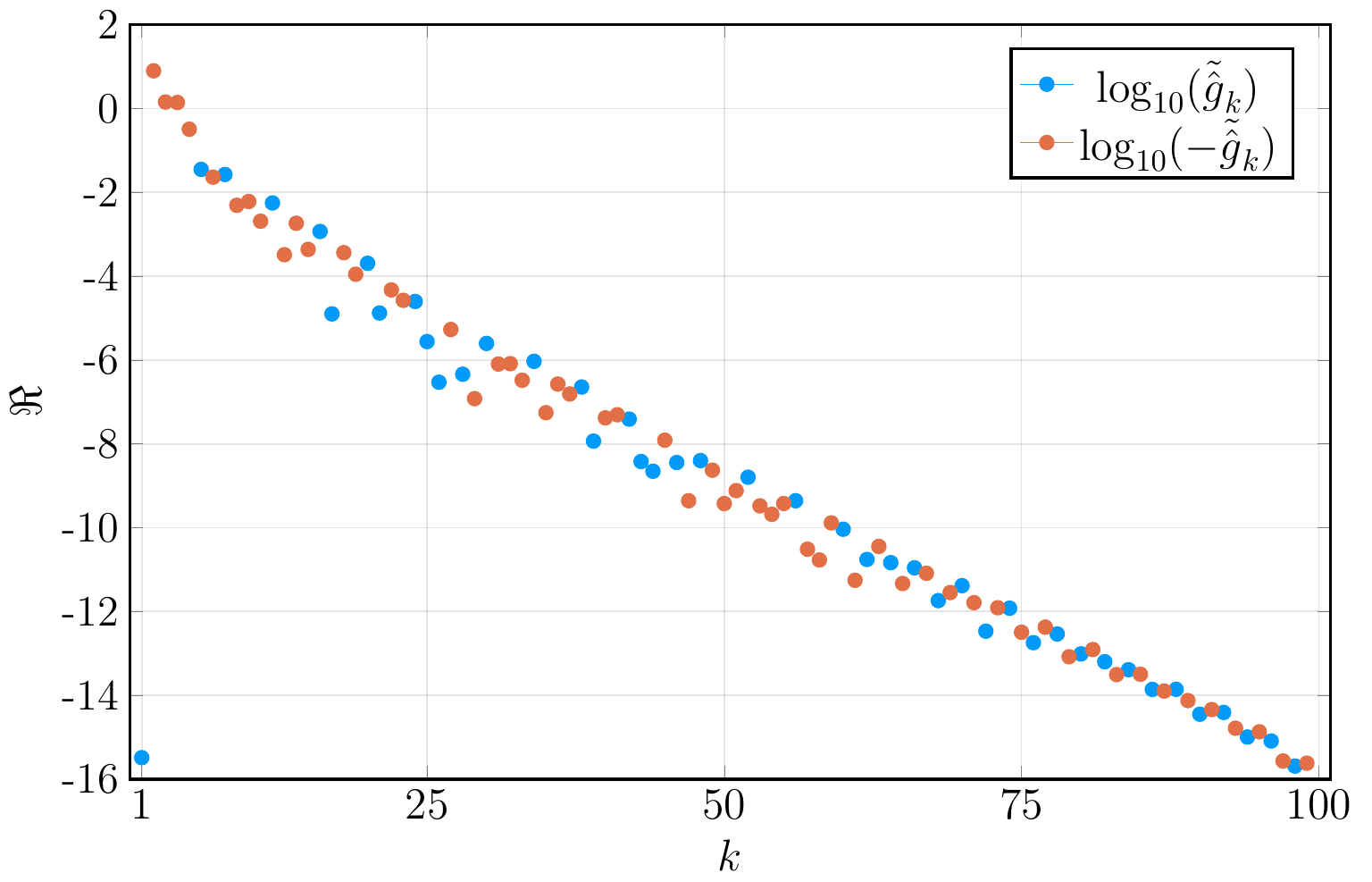}

\caption{[Example~\ref{exmp:8}: Symbol $f(\theta)=
 e^{3\mathbf{i}\theta}
-e^{2\mathbf{i}\theta}
+7e^{\mathbf{i}\theta}
+9e^{-1\mathbf{i}\theta}
-2e^{-2\mathbf{i}\theta}
+2e^{-3\mathbf{i}\theta}
-e^{-4\mathbf{i}\theta}$] Left: The approximated expansion functions $\tilde{c}_k(\theta_{j,n_0}), k=0,\ldots,\alpha$ for $n_0=100$ and $\alpha=4$. The approximation $\tilde{c}_0(\theta_{j,n_0})$ correspond well with $\lambda_j(T_{1000}(f))$ in the right panel of Figure~\ref{fig:exmp:4:symbols} (since $g$ is unknown). Right: The absolute value of the approximated first one houndred Fourier coefficients, $\log_{10}|\tilde{\hat{g}}_k|$.}
\label{fig:exmp:8:expansion}
\end{figure}
\begin{table}[!ht]
\centering
\caption{[Example~\ref{exmp:8}: Symbol $f(\theta)=
 e^{3\mathbf{i}\theta}
-e^{2\mathbf{i}\theta}
+7e^{\mathbf{i}\theta}
+9e^{-1\mathbf{i}\theta}
-2e^{-2\mathbf{i}\theta}
+2e^{-3\mathbf{i}\theta}
-e^{-4\mathbf{i}\theta}$] First ten computed ($\tilde{\hat{g}}_k$) Fourier coefficients of the unknown $g$. Approximations computed using $n_0=100$ and $\alpha=4$, and 512 bit precision.}
\label{tbl:exmp:8}
\begin{tabular}{rrrr}
\toprule
k&$\tilde{\hat{g}}_k$\\
\midrule
0& $-0.000000000000003$\\
1& $-7.931536795875190$\\
2& $-1.429849731406187$\\
3& $-1.393034471115375$\\
4& $-0.321121280053002$\\
5& $ 0.035288447846840$\\
6& $-0.023038821632295$\\
7& $ 0.026692519463291$\\
8& $-0.004916206049977$\\
9& $-0.006047350374789$\\
 \bottomrule
\end{tabular}
\end{table}
\end{exmp}

\section{Conclusions}
\label{sec:conclusions}
The working hypothesis in this article concerns the existence of an asymptotic expansion, such that there exists of a function $g$ describing the eigenvalue distribution of the Toeplitz matrices $T_n(f)$ generated by a symbol $f$.
We have shown numerically that we can recover an approximation of the function $g$. 
This is done by a matrix-less method described in Algorithm~\ref{algo:1}, which in principle can be modified so as to work without any information on $f$ or the way in which the eigenvalues of the smaller versions of $T_n(f)$ are computed. Algorithm~\ref{algo:1} can also be used to fast and accurately compute the eigenvalues of $T_n(f)$ for an arbitrarily large order $n$, as highlighted in \eqref{eq:describing:approximating:expansion}.
However, in this article we have focused only on using the obtained approximation of $g$ to find an approximation of its truncated Fourier series; and if $g$ is an RCTP, we have shown that it we are able to recover the original function $g$ analytically.
These approaches can be a valuable tool for the exploration of the spectrum of Toeplitz and Toeplitz-like  matrices previously not easily understood, because of high computational cost. 
For future research we propose the extension to complex-valued functions $g$ of the results presented herein, and also the study of matrices more general than $T_n(f)$. 
\section{Acknowledgments}
The author would like to thank Carlo Garoni and Stefano Serra-Capizzano for valuable insights and suggestions during the preparation of this work. The author is financed by Athens University of Economics and Business.
\bibliography{References}

\begin{thebibliography}{10}

\bibitem{ahmad171}
{\sc F.~Ahmad, E.~S. Al-Aidarous, D.~A. Alrehaili, S.-E. Ekstr\"{o}m, I.~Furci,
  and S.~Serra-Capizzano}, {\em Are the eigenvalues of preconditioned banded
  symmetric {T}oeplitz matrices known in almost closed form?}, Numerical
  Algorithms, 78 (2017), pp.~867--893.

\bibitem{barrera181}
{\sc M.~Barrera, A.~B\"{o}ttcher, S.~M. Grudsky, and E.~A. Maximenko}, {\em
  Eigenvalues of even very nice {T}oeplitz matrices can be unexpectedly
  erratic}, in The Diversity and Beauty of Applied Operator Theory, Springer
  International Publishing, 2018, pp.~51--77.

\bibitem{beam931}
{\sc R.~M. Beam and R.~F. Warming}, {\em The {A}symptotic {S}pectra of {B}anded
  {T}oeplitz and {Q}uasi-{T}oeplitz {M}atrices}, {SIAM} Journal on Scientific
  Computing, 14 (1993), pp.~971--1006.

\bibitem{bezanson171}
{\sc J.~Bezanson, A.~Edelman, S.~Karpinski, and V.~B. Shah}, {\em Julia: {A}
  {F}resh {A}pproach to {N}umerical {C}omputing}, {SIAM} Review, 59 (2017),
  pp.~65--98.

\bibitem{bogoya151}
{\sc J.~M. Bogoya, A.~B\"{o}ttcher, S.~M. Grudsky, and E.~A. Maximenko}, {\em
  Eigenvalues of {H}ermitian {T}oeplitz matrices with smooth simple-loop
  symbols}, Journal of Mathematical Analysis and Applications, 422 (2015),
  pp.~1308--1334.

\bibitem{bogoya171}
{\sc J.~M. Bogoya, S.~M. Grudsky, and E.~A. Maximenko}, {\em Eigenvalues of
  {H}ermitian {T}oeplitz {M}atrices {G}enerated by {S}imple-loop {S}ymbols with
  {R}elaxed {S}moothness}, in Large Truncated Toeplitz Matrices, Toeplitz
  Operators, and Related Topics, Springer International Publishing, 2017,
  pp.~179--212.

\bibitem{bottcher101}
{\sc A.~B\"{o}ttcher, S.~M. Grudsky, and E.~A. Maxsimenko}, {\em Inside the
  eigenvalues of certain {H}ermitian {T}oeplitz band matrices}, Journal of
  Computational and Applied Mathematics, 233 (2010), pp.~2245--2264.

\bibitem{bottcher991}
{\sc A.~B\"{o}ttcher and B.~Silbermann}, {\em Introduction to {L}arge
  {T}runcated {T}oeplitz {M}atrices}, Springer New York, 1999.

\bibitem{ekstrom191}
{\sc S.-E. Ekstr\"{o}m}, {\em Approximating the {P}erfect {S}ampling {G}rids
  for {C}omputing the {E}igenvalues of {T}oeplitz-like {M}atrices {U}sing the
  {S}pectral {S}ymbol}, 2019.
\newblock arXiv:1901.06917 (submitted).

\bibitem{ekstrom185}
{\sc S.-E. Ekstr\"{o}m, I.~Furci, C.~Garoni, C.~Manni, S.~Serra-Capizzano, and
  H.~Speleers}, {\em Are the eigenvalues of the {B}-spline isogeometric
  analysis approximation of {$-\Delta u=\lambda u$} known in almost closed
  form?}, Numerical Linear Algebra with Applications, 25 (2018), p.~e2198.

\bibitem{ekstrom181}
{\sc S.-E. Ekstr\"{o}m, I.~Furci, and S.~Serra-Capizzano}, {\em Exact formulae
  and matrix-less eigensolvers for block banded symmetric {T}oeplitz matrices},
  {BIT} Numerical Mathematics, 58 (2018), pp.~937--968.

\bibitem{ekstrom183}
{\sc S.-E. Ekstr\"{o}m and C.~Garoni}, {\em A matrix-less and parallel
  interpolation--extrapolation algorithm for computing the eigenvalues of
  preconditioned banded symmetric {T}oeplitz matrices}, Numerical Algorithms,
  (2018).
\newblock https://doi.org/10.1007/s11075-018-0508-0 (in press).

\bibitem{ekstrom171}
{\sc S.-E. Ekstr\"{o}m, C.~Garoni, and S.~Serra-Capizzano}, {\em Are the
  {E}igenvalues of {B}anded {S}ymmetric {T}oeplitz {M}atrices {K}nown in
  {A}lmost {C}losed {F}orm?}, Experimental Mathematics, 27 (2017),
  pp.~478--487.

\bibitem{ekstrom184}
{\sc S.-E. Ekström}, {\em Matrix-{L}ess {M}ethods for {C}omputing
  {E}igenvalues of {L}arge {S}tructured {M}atrices}, PhD thesis, Uppsala
  University, Uppsala: Acta Universitatis Upsaliensis, 2018.

\bibitem{garoni171}
{\sc C.~Garoni and S.~Serra-Capizzano}, {\em Generalized {L}ocally {T}oeplitz
  {S}equences: {T}heory and {A}pplications (Volume 1)}, Springer International
  Publishing, 2017.

\bibitem{noack191}
{\sc A.~Noack}, {\em Generic{L}inear{A}lgebra.jl}.
\newblock https://github.com/JuliaLinearAlgebra/GenericLinearAlgebra.jl.

\bibitem{parter621}
{\sc S.~V. Parter and J.~Youngs}, {\em The symmetrization of matrices by
  diagonal matrices}, Journal of Mathematical Analysis and Applications, 4
  (1962), pp.~102--110.

\bibitem{reichel921}
{\sc L.~Reichel and L.~N. Trefethen}, {\em Eigenvalues and pseudo-eigenvalues
  of {T}oeplitz matrices}, Linear Algebra and its Applications, 162-164 (1992),
  pp.~153--185.

\bibitem{shapiro171}
{\sc B.~Shapiro and F.~{\v{S}}tampach}, {\em Non-{S}elf-{A}djoint {T}oeplitz
  {M}atrices {W}hose {P}rincipal {S}ubmatrices {H}ave {R}eal {S}pectrum},
  Constructive Approximation,  (2017).
\newblock https://doi.org/10.1007/s00365-017-9408-0 (in press).

\bibitem{tilli991}
{\sc P.~Tilli}, {\em Some results on complex {T}oeplitz eigenvalues}, Journal
  of Mathematical Analysis and Applications, 239 (1999), pp.~390--401.

\bibitem{trefethen051}
{\sc L.~N. Trefethen and M.~Embree}, {\em Spectra and pseudospectra: the
  behavior of nonnormal matrices and operators}, Princeton University Press,
  Princeton, N.J, 2005.

\end{thebibliography}
\bibliographystyle{siam}
\end{document}